\documentclass[reqno]{amsart}

\usepackage{graphicx,helvet}
\usepackage{epstopdf}
\usepackage[percent]{overpic}
\usepackage{xcolor}
\usepackage{multirow}

\usepackage{epsfig}
\usepackage{cite}
\usepackage[english]{babel}
\usepackage[utf8]{inputenc}

\allowdisplaybreaks

\newcommand{\e}{\mathrm{e}}

\newcommand{\eps}{\varepsilon}
\renewcommand{\epsilon}{\varepsilon}

\newcommand{\R}{\mathbb{R}}

\newcounter{dctr}[section]
\numberwithin{equation}{section}

\newtheorem{remark}{Remark}[section]

\def\bcdot{{*}}

\allowdisplaybreaks
\numberwithin{equation}{section}

\title[ FOWENO-APM methods]{Lax Wendroff  approximate Taylor methods with fast and optimized weighted essentially non-oscillatory  reconstructions}
\author[Carrillo]{H. Carrillo$^{\mathrm{a}}$}
\author[Parés]{C. Parés$^{\mathrm{b}}$}
\author[Zorío]{D. Zorío$^{\mathrm{c}}$}

\usepackage[normalem]{ulem}

\begin{document}

\begin{abstract}
 The goal of this work is to introduce new families of shock-capturing high-order numerical methods for systems of conservation laws that combine Fast WENO (FWENO) and Optimal WENO (OWENO) reconstructions with Approximate Taylor methods for the time discretization. FWENO reconstructions are based on smoothness indicators that require a lower number of calculations than the standard ones. OWENO reconstructions are based on a definition of the nonlinear weights that allows one to unconditionally attain the optimal order of accuracy regardless of the order of critical points. Approximate Taylor methods update the numerical solutions by using a Taylor expansion in time in which, instead of using the Cauchy-Kovalevskaya procedure, the time derivatives are computed by combining spatial and temporal numerical differentiation with Taylor expansions in a recursive way. These new methods are compared between them and against methods based on standard WENO implementations and/or SSP-RK time discretization. A number of test cases are considered ranging from scalar linear 1d problems to nonlinear systems of conservation laws in 2d.
\end{abstract} 

\date{\today}


\thanks{$^{\mathrm{A}}$Departamento de Matemática Aplicada, 
Universidad de Málaga, Avda. Cervantes, 2. 29071 Málaga,  
 Spain.  E-Mail: 
   {\tt hcarrillo@uma.es}}

 \thanks{$^{\mathrm{B}}$Departamento de Matemática Aplicada, 
Universidad de Málaga, Avda. Cervantes, 2. 29071 Málaga,  
 Spain.  E-Mail: 
   {\tt pares@uma.es}}

\thanks{$^{\mathrm{C}}$CI$^{\mathrm{2}}$MA,   Universidad de
Concepci\'{o}n, Casilla 160-C, Concepci\'{o}n, Chile.  E-Mail:
  {\tt dzorio@ci2ma.udec.cl}}

\maketitle

\section{Introduction}\label{sec1}

Weighted Essentially Non-Oscillatory (WENO) reconstructions (see \cite{LiuWENO}, \cite{JiangShu96}) and
Strong Stability Preserving Runge-Kutta time discretizations (see \cite{Gottlieb1998}, \cite{Gottlieb2011}) have become common ingredients of high-resolution schemes for the numerical solution of hyperbolic conservation laws:
\begin{equation} \label{cons_law_0}
u_t + f(u)_x = 0 , \quad \quad u(x,0)=u_0(x),  \quad -\infty <x < \infty.
\end{equation}
Here $u:\R \times \R \rightarrow \R^m$ is an $m$-dimensional vector of conserved quantities. 

WENO methods present a high order of accuracy in smooth zones and avoid oscillatory behaviours close to discontinuities through the construction of non-linear weights based on some smooth indicators. Many variants of the original WENO reconstruction have been introduced since then. For instance, in FWENO methods introduced in \cite{ZBBM2019}, new smoothness indicators have been proposed that require a lower number of calculations than the ones proposed by Jiang and Shu.  

On the other hand, the expression of the weights in the original WENO method leads to an undesired loss of accuracy near critical points. Different variants have been introduced to deal with this difficulty: see \cite{Arandiga2011}, \cite{Arandiga2014}, \cite{HENRICK2005542}, \cite{YamaleevCarpenter2009}.
To the best of our knowledge the only approach that allows one to unconditionally attain the optimal order of accuracy regardless of the order of critical points is, for third order reconstructions, the OWENO3 method introduced in  \cite{BBMZO3} and, for reconstructions of order bigger than 3, the OWENO methods presented in \cite{BBMZSINUM}. In this latter reference, the Jiang-Shu smoothness indicators are used to define the weights (for third order methods these indicators coincide with those of FWENO methods). 
In this work, the following WENO reconstructions will be used:
\begin{itemize}
    \item OWENO3 method for third order reconsructions;
    \item WENO methods based on the expression of the OWENO weights and the smoothness indicators of FWENO,
    so that they are both fast and optimal.
\end{itemize}
For shortness, we will refer to these methods as FOWENO reconstructions. 

Concerning the time stepping, an alternative to SSP-RK methods is given by methods that use Taylor expansions in time to update the numerical solution \begin{equation}\label{atm_0}
u_i^{n+1} = u_i^n +\sum_{k= 1}^{m} \frac{\Delta t ^k}{k! }\ u^{k}_i + \mathcal{O}\left(\Delta t^{m+1}\right).
\end{equation}
where $\{ x_i \}$ are the nodes of a uniform mesh of step $\Delta x$; $u_i^n$ is an approximation of the point value of the solution at $x_i$ at the time $n \Delta t$, where $\Delta t$ is the time step; and $u^k_i$ is an approximation of the $k$-order time derivative of $u$ at $x_i$ at time $n \Delta t$. Although the values of $u_i^{k}$  can be approximated using the  Cauchy - Kovalevskaya (CK) procedure, it is well-known that, for nonlinear problems, this approach may be impractical from the computational point of view (symbolic calculus, tensor matrix, excessive computations...) In the context of ADER methods introduced by Toro and collaborators (see \cite{Ader2001},  \cite{TitarevToro2002},   \cite{Schwartzkopff2002}), this difficulty have been circumvented   by replacing the CK procedure by local space-time problems that are solved with a Galerkin method: see \cite{DumbserToro2008}, \cite{PNPM}. 

We follow here the strategy introduced in \cite{ZBM2017} to avoid the CK procedure. It is based on the equalities  
\begin{equation}\label{u_t^k}
\partial_t^k u = -\partial_x\partial_t^{k-1}f(u).
\end{equation}
that can be easily derived from the equation, if the solutions are assumed to be smooth enough. Numerical approximations of the derivatives appearing at the right-hand side are computed by combining numerical differentiation formulas in space and time with Taylor expansions in a recursive way. 

The so-called Lax-Wendroff approximate Taylor (LAT) methods introduced in \cite{ZBM2017} do not generalize the standard Lax-Wendroff methods for linear systems: if, for instance, a LAT method that updates the numerical solution using  \eqref{atm_0} with 
$m = 2$ and uses 3-point centered formulas to approximate the derivatives is applied to  \eqref{cons_law_0} with $f(u) = au$, the  numerical scheme obtained is
\begin{equation}\label{5plwm}
u_i^{n+1} = u_i^n - \frac{a \Delta t}{2\Delta x}(u_{i+1}^n - u_{i-1}^n)  - \frac{a^2\Delta t^2}{8\Delta x^2}(u_{i+2}^n - 2 u_i^n  + u_{i-2}^n),
\end{equation}

which is different from the standard Lax-Wendroff method and whose stability properties  are worse (see \cite{LeVeque2007book}). Compact Approximate Taylor (CAT) methods were designed in \cite{CP2019} as a variant of these methods that properly generalize  the Lax-Wendroff methods for linear systems.  

Although both LAT and CAT strategies have been combined previously with standard WENO reconstructions as equipment for attaining shock-capturing properties (see \cite{ZBM2017} and \cite{CP2019}), they have been never combined with FOWENO reconstructions: the goal of this work is to introduce new families of high-order numerical methods using FOWENO reconstructions and Approximate Taylor methods. These methods will be  compared between them and against standard WENO implementations in a number of test cases ranging from scalar linear 1d problems to nonlinear systems of conservation laws in 2d.

The paper is organized as follows. In section \ref{sec2},  the LAT and CAT strategies to derive approximate Taylor methods are briefly recalled. In section \ref{foweno}, the new FOWENO reconstructions are described in detail. In section \ref{sec4}, the ingredients already described in section \ref{sec2} and \ref{foweno} are combined to construct FOWENO-APT methods. Section \ref{sec5} focuses on the numerical experiments:
methods based on WENO or FOWENO reconstructions combined with CAT, LAT or SSPRK are applied to the 1d linear transport equation, Burgers equation, and the 1d and 2d Euler equations of gas dynamics. The quality of the solutions and the CPU run-time are compared and discussed. Finally, in section \ref{sec6} some conclusions are given.

\section{Approximate Taylor Methods} \label{sec2}

Approximate Taylor methods are based on  a Taylor expansion in time (\ref{atm_0}) to update the numerical solutions in which time
derivatives are computed by using the equalities \eqref{u_t^k}. For  the  sake  of  simplicity,  the  methods  will  be  only  described  for  the one-dimensional scalar case.

\subsection{Lax-Wendroff Approximate Taylor Methods }\label{LAT} 

In Lax-Wendroff Approximate Taylor(LAT) methods,  the time derivatives $\partial_t^{k}u$ are approximated
by applying a first order numerical differentiation formula in space to some approximations 
\begin{equation}\label{f(k-1)_i}
\tilde{f}^{(k-1)}_{i}  \cong \partial_t^{k-1}f(u)(x_{i}, t_n)
\end{equation}
that will be computed by using recursively Taylor expansions in time. 

LAT methods are based on centered $(2p +1)$-point numerical differentiation formulas
\begin{equation} \label{F}
f^{(k)}(x_i)   \simeq   D^k_{p,i}(f, \Delta x) =  \frac{1}{\Delta x^k} \sum_{j=-p}^{p} \delta^k_{p,j} f(x_{i+j}).
\end{equation}
The following  notation 
\begin{equation}
 D^k_{p,i}(f_\bcdot, \Delta x)   =   \frac{1}{\Delta x^k} \sum_{j=-p}^{p} \delta^k_{p,j} f_{i+j},
\end{equation}
will be used to indicate that the formula is applied to some approximations $f_i$ of $f$ and not to its exact point values $f(x_i)$. 
In cases where there are two or more indexes, the symbol $\bcdot$ will be used to indicate with respect to which the differentiation is applied. For instance:
\begin{eqnarray*}
\partial^k_x u(x_i , t_n) & \simeq & D^{k}_{p,i} (u_\bcdot^n, \Delta x) =  \frac{1}{\Delta x^k} \sum_{j = -p}^p 
\delta^{k}_{p,j} u^n_{i+j}, \\
\partial^k_t u(x_i, t_n) &\simeq& D^{k}_{p,n} (u_i^\bcdot, \Delta t) =  \frac{1}{\Delta t^k} \sum_{r = -p}^p \delta^{k}_{p,r} u^{n+r}_{i}.
\end{eqnarray*}

Once the approximations \eqref{f(k-1)_i} have been computed, the time derivatives of the solution are approximated by:
$$
\partial_t^k u(x_i, t_n) \cong  \tilde u^{(k)}_i =  - D^1_{p,i} (\tilde{f}^{(k-1)}_\bcdot, \Delta x) = - \frac{1}{\Delta x}\ \sum_{j=-p}^p \delta^1_{p,j}  \tilde{f}^{(k-1)}_{i+j}.
$$

A recursive procedure is followed to compute the approximation of the time derivatives: once $u^l_i$, $l=0, \dots, k$ have been computed, 
a Taylor expansion of degree $k$ is used to compute approximations $\tilde{f}^{k-1,n+r}_{i}$ of  $f(u(x_i, (n+r)\Delta t)$, $r = -p, \dots, p$; the centered differentiation formula is then used to obtain $\tilde{f}^{(k-1)}_{i}$; and, finally, the first order derivative in space is applied to $\tilde{f}^{(k-1)}_{i + j}$, $j = -p, \dots, p$  to compute $u_i^{k+1}$.  Once all the time derivatives are approximated,  \eqref{atm_0} is used to update the numerical solutions. 

The procedure can be summarized as follows:

\begin{enumerate}

\item Define
$$
\tilde{f}^{(0)}_i = f(u^n_{i}).
$$

\item Compute
 \begin{equation}\label{lat_ut}
 \tilde{u}  ^{(1)}_{i} = - D^1_{p,i}(\tilde{f}^{(0)}_\bcdot, \Delta x).
 \end{equation}

\item  For $k = 2, \dots, m$:

\begin{enumerate}

\item Compute  
$$
\tilde{f}^{k-1,n+r}_{i} = f \left(  u^n_{i} + \sum_{l=1}^{k-1} \frac{(r \Delta t)^l}{l!} \tilde{u}^{(l)}_{i} \right), \quad 
r = -p, \dots, p.
$$

\item Compute 
\begin{equation}\label{derf(k-1)}
\tilde{f}^{(k-1)}_{i} =  D^{k-1}_p(\tilde{f}^{k-1, \bcdot}_i, \Delta t).
\end{equation}

\item Compute 
\begin{equation}\label{deru(k)}
    \tilde{u}  ^{(k)}_{i} = - D^1_{p,i}(\tilde{f}^{(k-1)}_\bcdot, \Delta x).
\end{equation}

\end{enumerate}

\item Update the solution by (\ref{atm_0}).

\end{enumerate}

The order of the method is $\min(m, 2p)$.

 \begin{remark}
 Although, for the sake of clarity, $m$ and $p$ have been considered as two arbitrary positive integers in the presentation of LAT methods, in \cite{ZBM2017} 
 $m$ is an odd number (since the method is combined with WENO reconstructions) and $p$ is chosen adequately to obtain order $m$. More precisely, in formulas \eqref{deru(k)},
 $$
 p = \left\lceil  \frac{m+1-k}{2}  \right\rceil,
 $$
 where $\lceil \cdot \rceil$ is the ceiling function, and in formulas \eqref{derf(k-1)}
 $$
 p = \frac{m- 1}{2}.
 $$
 \end{remark}







LAT methods can be written in conservative form. To see this, let us introduce the family of interpolatory numerical differentiation formulas
\begin{equation}\label{upwF}
f^{(k)}(x_i + q \Delta x) \simeq A^{k,q}_{p,i}(f, \Delta x) = \frac{1}{\Delta x^k} \sum_{j = -p + 1}^p \gamma^{k,q}_{p,j} f(x_{i+j}),
\end{equation}
 that approximates the $k$-th derivative of a function at the point $x_i + q \Delta x$ using its values at the $2p$  points $x_{i-p+1}, \dots, x_{i+p}$. The symbol $\bcdot$ will be used again to indicate whit respect to which index the differentiation is performed. 

The following relation holds (see \cite{CP2019}):
\begin{equation}\label{F2}
D^k_{p,i}(f, \Delta x)   =  \frac{1}{\Delta x} \left(A^{k-1, 1/2}_{p,i}(f, \Delta x) -  A^{k-1, 1/2}_{p,i-1}(f, \Delta x)\right).
\end{equation}
Using this equality with $k = 1$, LAT methods can be written in the form
  \begin{equation}\label{cons}
u_i^{n+1} = u_i^n + \frac{\Delta t}{\Delta x}\left( F^p_{i-1/2} - F^p_{i+1/2}\right),
\end{equation}
where
\begin{equation}\label{cat1}
F^p_{i+1/2}  = \sum_{k=1}^{m} \frac{\Delta t^{k-1}}{k!}A^{0, 1/2}_{p,i}(\tilde{f}_{i,\bcdot}^{(k-1)}, \Delta x).  
\end{equation}

\subsection{Compact Approximate Taylor methods}

CAT methods are based on  the conservative expression \eqref{cons}-\eqref{cat1}, with the difference that now only the values
\begin{equation} \label{fluxstencil}
u^n_{i -p +1}, \dots, u^n_{i+ p}, 
\end{equation}
are used to compute
the numerical flux $F_{i+1/2}$, so that a centered $(2p + 1)$-point stencil is used to compute $u_i^{n+1}$. 
The numerical flux is thus computed as follows:
\begin{equation}\label{cat2}
F^p_{i+1/2}  = \sum_{k=1}^{m} \frac{\Delta t^{k-1}}{k!}A^{0, 1/2}_{p,0}(\tilde{f}_{i,\bcdot}^{(k-1)}, \Delta x).  
\end{equation}
where 
\begin{equation}\label{f(k-1)_ij}
\tilde{f}^{(k-1)}_{i,j}  \cong \partial_t^{k-1}f(u)(x_{i+j}, t_n), \quad j=-p+1,\dots, p
\end{equation}
are \textit{local} approximations of the time derivatives of the flux.  By \textit{local} we mean that these approximations depend on the stencil, i.e.
$$
i_1 + j_1 = i_2 + j_2  \not \Rightarrow \tilde{f}^{(k-1)}_{i_1,j_1} = \tilde{f}^{(k-1)}_{i_2,j_2}.
$$
Local approximations of the time derivatives of the solution
\begin{equation*}
\tilde{u}  ^{(k)}_{i,j} \cong \partial_t^{(k)}u(x_{i+j}, t_n), \quad j=-p+1,\dots, p
\end{equation*}
are obtained then by using the non-centered differentiation formulas
\begin{equation*}
\tilde{u} ^{(k)}_{i,j} = - A^{1,j}_{p,0}(\tilde{f}^{(k-1)}_{i, \bcdot}, \Delta x)
= - \frac{1}{\Delta x} \sum_{r=-p +1}^p \gamma^{1,j}_{p,r} \tilde{f}^{(k-1)}_{i, r}.
\end{equation*}

Like in LAT methods, these local approximations of the time derivatives are recursively used to compute approximations of the flux forward and backward in time using Taylor expansions in a recursive way. 

Given $i$, the procedure to compute $F^p_{i+1/2}$ is as follows:

\begin{enumerate}

\item  {Define}
$$
\tilde{f}^{(0)}_{i,j}=f(u^n_{i+j}), \quad  j = -p+1, \dots, p.
$$

\item  {For $k = 2 \dots m$:}

\begin{enumerate}

\item Compute
\begin{equation*}
 \tilde{u}^{(k-1)}_{i,j} = - A^{1,j}_{p,0}(\tilde{f}^{(k-2)}_{i,\bcdot}, \Delta x). 
\end{equation*}

\item Compute 
$$
\tilde{f}^{k-1,n+r}_{i,j} = f \left(  u^n_{i+j} + \sum_{l=1}^{k-1} \frac{(r \Delta t)^l}{l!} \tilde{u}^{(l)}_{i,j} \right), \quad  j, r = -p+1, \dots, p.
$$

\item Compute
$$
\tilde{f}^{(k-1)}_{i,j} =   A^{k-1,0}_{p,n}( \tilde{f}^{k-1, \bcdot}_{i,j}, \Delta t),\quad  j = -p+1, \dots, p.
$$

\end{enumerate}

\item Compute $F^p_{i+1/2}$ by \eqref{cat2}

\end{enumerate}

\noindent Once the numerical fluxes have been computed, the numerical solution is updated by using \eqref{cons}.

In \cite{CP2019} it has been shown that:

\begin{itemize}
    \item The order of the method is $\min(m, 2p)$ so that the optimal choice is $m = 2p$: the corresponding numerical method will be represented by CAT$2p$ in the sequel.
    \item CAT$2p$ reduces to the standard Lax-Wendroff method for linear problems.
    \item CAT$2p$ is linearly stable under the standard CFL-1 condition.
    
\end{itemize}

\noindent The extension of LAT and CAT methods to systems is straightforward by applying the schemes component by component. 
The extension to multiple dimensions using Cartesian grids can be done through the methods of lines. For a 2D problem, CAT uses a rectangular stencil of $p^2$ points centered in a point 
$(x_{i + 1/2}, y_{j + 1/2})$ to compute the horizontal component of the  numerical flux at $(x_{i+ 1/2}, y_j) $ and the vertical component at $(x_i, y_{j + 1/2})$ on the basis of local approximations of  the time derivatives and applications of Taylor expansions.

\section{Fast and optimal WENO reconstructions}\label{foweno}
Approximate Taylor methods produce spurious oscillations near discontinuities due to the Gibbs phenomenon. In order to get rid of these oscillations, WENO reconstructions will be used to compute the first order derivatives in time. 

Given the point values of a function $f$ at a stencil of $2p+1$ points:
\begin{align*}
S_i=\{f_{i-p},\ldots,f_{i+p}\},
\end{align*}
where 
$f_j = f(x_j)$, WENO operators provide a reconstruction of $f$ at 
$$x_{i+1/2} = x_i + \frac{h}{2}, $$
where $h$ is the step of the mesh (assumed to be constant).  This reconstruction is based on the Lagrange interpolation polynomials  $p_{s}(x)$, $0\leq s\leq p$ that interpolates the point values at $p+1$ sub-stencils
\begin{align*}
S_{p,s} = \{ f_{i-p+s}, \dots, f_{i+s} \}, \quad s = 0, \dots, p.
\end{align*}
More precisely, the WENO strategy consists in defining the reconstruction as a convex combination
$$
q(x_{i+1/2}) = \sum_{s=0}^p w_s p_s(x_{i+1/2}),
$$
where the weights $w_0, \dots, w_p $ satisfy $w_s \cong c_s$ on smooth zones, where $c_0, \dots, c_p$ are the linear ideal weights satisfying 
$$
P(x_{i+1/2}) = \sum_{s=0}^p c_s p_s(x_{i+1/2}),
$$
where $P(x)$ is the polynomial that interpolates all the point values of the stencil $S_i$. The weights $w_i$ are function of some smoothness indicators. In FWENO methods introduced in \cite{ZBBM2019}, the following smoothness indicators have been proposed
\begin{equation}\label{smind} 
I_{s}:=\sum_{j=1}^{p}(f_{-p+i+s}-f_{-p-1+i+s})^2,\quad0\leq s\leq p,
\end{equation} 
that require a lower number of calculations than the smoothness indicators by Jiang and Shu (see \cite{JiangShu96}).  

On the other hand, the expression of the weights in the original WENO method leads to an undesired loss of accuracy near critical points. To the best of our knowledge the only approach that allows to unconditionally attain the optimal order of accuracy regardless of the order of critical points is, for third order reconstructions, the OWENO3 method introduced in  \cite{BBMZO3} and, for reconstructions of order bigger than 3, the OWENO methods presented in \cite{BBMZSINUM}. In this latter reference, the Jiang-Shu smoothness indicators are used to define the weights (for third order methods these indicators coincide with \eqref{smind}).

Let us summarize here the expression of FOWENO methods 
(see  \cite{BBMZO3} and \cite{BBMZSINUM} for the accuracy analysis). The expression of FOWENO3, (i.e. OWENO$3$) is the following:

Given $i$ and $\varepsilon>0$,

\smallskip

\begin{enumerate}

\item Increase the dependence data stencil 
\begin{equation}
\bar{S}=\{f_{i-1},f_i,f_{i+1},f_{i+2}\},
\end{equation}
with $f_i=f(x_i)$.

\item Compute the corresponding interpolating polynomials evaluated at
  $x_{i+1/2}$, which, both in case of reconstructions from point values
  and from cell averages, are given by
  \begin{align} 
  p_0(x_{i+1/2}) = -\frac{1}{2} f_{i-1} + \frac{3}{2} f_i,\quad 
  p_1(x_{i+1/2}) = \frac{1}{2} f_i + \frac{1}{2} f_{i+1}.
\end{align} 
  
\item Compute the corresponding Jiang-Shu smoothness indicators $I_0$,
  $I_1$ and~$I_2$ (including the one considering the rightmost node)
  by
  \begin{align}
  I_0 = (f_i-f_{i-1})^2,\quad 
  I_1 = (f_{i+1}-f_i)^2,\quad
  I_2 = (f_{i+2}-f_{i+1})^2.
\end{align}

\item Compute the preliminary weights~$\tilde{\omega}_0$
  and~$\tilde{\omega}_1$:
  \begin{align} \label{eq:tw}  
    \tilde{\omega}_s := \frac{I_s+\varepsilon}{I_0+I_1+2\varepsilon},\quad 
    s=0,1
  \end{align} 

\item Define $\tau$ by
  \begin{align} \label{eq:taudef} 
    \tau:=dI,\quad d:=(-f_{i-1}+3f_i-3f_{i+1}+f_{i+2})^2,\quad I:=I_0+I_1+I_2. 
  \end{align} 
\item Compute the corrector weight~$\omega$:
  \begin{align} \label{eq:cw} 
    \omega=\frac{J}{J+\tau+\varepsilon}, \quad 
    \text{with $J=I_0(I_1+I_2)+(I_0+I_1)I_2$.} 
  \end{align} 
\item Compute the corrected weights $\omega_0$ and $\omega_1$:
  \begin{align} \label{correctedweights}  
&   \omega_0 := \omega c_0+(1-\omega)\tilde{\omega}_0,\quad 
  \omega_1 :=  \omega c_1+(1-\omega)\tilde{\omega}_1,
 \end{align}
 where $c_0$, $c_1$ are the ideal linear weights.
\item Obtain the OWENO reconstruction at $x_{i+1/2}$:
  $$q(x_{i+1/2})=\omega_0p_0(x_{i+1/2})+\omega_1p_1(x_{i+1/2}).$$
\end{enumerate}

Unlike FOWENO3, FOWENO(2p+1)  reconstructions for $p\geq2$ do not require to increase artificially the stencil. Their expression, combined with the smoothness indicators \eqref{smind}  can be summarized as follows:

\smallskip 
 
Given $i$, the stencil $S_i$ and $\varepsilon>0$.

\begin{enumerate}
\item Compute the interpolating polynomials  $p_{j}$, $j = 0\leq j\leq p,$

\item Compute the fast smoothness indicators  \eqref{smind}.

\item {Compute the discriminant}
  
  \begin{equation*}
      D_p =|B_p-4A_p C_p|,
  \end{equation*}
with 
  \begin{align} \label{difnodiv2}
  A_p= \,\frac{1}{2} \sum_{j=-p}^{p} \delta^{2p}_{p,j}f_{i+j}, \quad
   B_p= \,\sum_{j=-p}^{p} \delta^{2p-1}_{p,j}f_{i+j},
  \quad
  C_p= \, \sum_{j=-p}^{p} \delta^{2p-2}_{p,j}f_{i+j}.
  \end{align}
for  $j=-p,\dots,p$.

\item Obtain the squared undivided difference of order $2p$:
  \begin{align} \label{difnodiv}
  \tau_p =(2A_p)^2.
  \end{align}

\item Compute 
\begin{equation*}
    d_p:=\frac{\tau_p^{a_1}D_p^{a_1}}{\tau_p^{a_1}+D_p^{a_1}+\eps}
\end{equation*}
  for some $a_1$ chosen by the user such that $a_1\geq1$, as done in \cite{BBMZSINUM}.
  
\item Compute  
  \begin{align} \label{alphayc} 
  \alpha_{s}=c_{s}\biggl(1+\frac{d_p}{I_{s}^{a_1}+\varepsilon}\biggr)^{a_2},\quad0\leq s\leq p, 
  \end{align} 
  where  $c_{s}$ are the ideal linear weights.  $a_2$ is
  chosen by the user such that $a_2\geq\frac{p+1}{2a_1}$, which is a
  sufficient condition to attain the optimal $(p+1)$-th accuracy near
  discontinuities \cite{ZBBM2019}.
  
\item Generate the FOWENO weights:
  \begin{align} \label{omegari} 
  \omega_{s}=\frac{\alpha_{s}}{\alpha_{0} + \dots+ \alpha_{p}}, \quad 
  s = 0, \dots, p. \end{align} 
\item Obtain the reconstruction at $x_{i+1/2}$:
  \begin{align} \label{qrx} 
  q_p(x_{i+1/2})=\sum_{s=0}^{p}\omega_{s}p_{s}(x_{i+1/2}). \end{align} 
\end{enumerate}

Combining the results obtained in \cite{ZBBM2019} and \cite{BBMZSINUM}
it follows that this method attains the optimal order regardless of
the order of the critical point, without having to artificially tune
$\varepsilon$.  

\section{FOWENO-ATM}\label{sec4}

With the FOWENO spatial reconstructions already defined, we incorporate them in the  Approximate Taylor methods to avoid the appearance of oscillations near the discontinuities or shocks,  substituting the first derivative in time of the Taylor expansion by those reconstructions.  More precisely, in LAT methods of Section (\ref{lat_ut}) is replaced by:

\begin{equation}\label{ut_foweno}
\tilde{u}  ^{(1)}_{t,i} = - \frac{\hat f_{i +  1/2} - \hat f_{i -1/2}}{\Delta x}.
\end{equation}
where $\hat f_{i + 1/2}$ denotes the $(2p+1)$-th order FOWENO flux splitting  reconstructions at $x_{i + 1/2}$.  In CAT methods,  \eqref{cat2} is replaced by:
\begin{equation}\label{catweno2}
F^p_{i+1/2}  = \hat f_{i +  1/2} + \sum_{k=2}^{m} \frac{\Delta t^{k-1}}{k!}A^{0, 1/2}_{p,0}(\tilde{f}_{i,\bcdot}^{(k-1)}, \Delta x).  
\end{equation}

FOWENO reconstructions are computed in conserved variables using the procedure described in  \cite{Shu1989}, so that their extension to systems is straightforward. 

\section{Numerical experiments} \label{sec5}

In order to simplify the notation and save space for the labels, from now on the following abbreviations will be used for the different numerical methods to be compared:

\begin{table}[htbp]
\begin{center}
\resizebox{9cm}{!}{
\begin{tabular}{|l||l|}
\hline 
 Abbreviation  & Numerical method \\
\hline \hline
W$q$R$s$ & WENO$q$ with SSPRK$s$ \\
\hline 
 W$q$C$s$ &  WENO$q$ with CAT$s$\\
\hline  
W$q$L$s$ &   WENO$q$ with LAT$s$\\
\hline  
FOW$q$R$s$ &  FOWENO$q$ with SSPRK$s$\\
\hline  
FOW$q$C$s$ &   FOWENO$q$ with CAT$s$\\
\hline  
FOW$q$L$s$ &   FOWENO$q$ with LAT$s$\\
\hline  
\end{tabular}
}
\vspace{2mm}
\label{metodos}
\end{center}
\end{table}

\noindent Here, SSPRK denotes the well-known Strong Stability Preserving Runge-Kutta methods \cite{Gottlieb2011}, $q$ is the accuracy order of the spatial WENO reconstructions and $s$ is the order of accuracy of the time discretization. We present some numerical experiments using FOWENO and the traditional WENO \cite{Shu1989} reconstructions combined with CAT$\{2,4,6\}$, LAT$\{3,5,7\}$  and SSPRK$\{3,4\}$  over some classical 1D scalar conservation laws (linear transport and Burgers equations) and 1D and 2D systems (Euler equations of gas dynamics).

\subsection{Scalar conservation laws} 
Let us consider first the one-dimensional scalar conservation law:
\begin{equation} \label{cons_law}
u_t + f(u)_x = 0.
\end{equation}
\subsubsection{Test 1: Linear transport equation}
We consider (\ref{cons_law}) with linear flux function $f(u)= au$ in the spatial interval $x\in[0,2]$ with initial condition:

\begin{equation} \label{initest1} 
u(x,0)= \left\{ \begin{array}{l r c l}
e^{-1200(x-1/3)^2}     & 0    &\leq x < &2/3, \\
6(x-2/3)               & 2/3  &\leq x < &5/6, \\
-6(x-1)                & 5/6  &\leq x < &1, \\
1                      & 7/6  &\leq x \leq& 4/3, \\
\sqrt{1-100(x-5/3)^2}  & 3/4  &< x \leq &2. \\
\end{array}\right.
\end{equation}
 
\begin{figure}[!ht]
   \small
	\setlength{\unitlength}{1mm}
	\centering	
	\includegraphics[width=\textwidth, height=9cm]{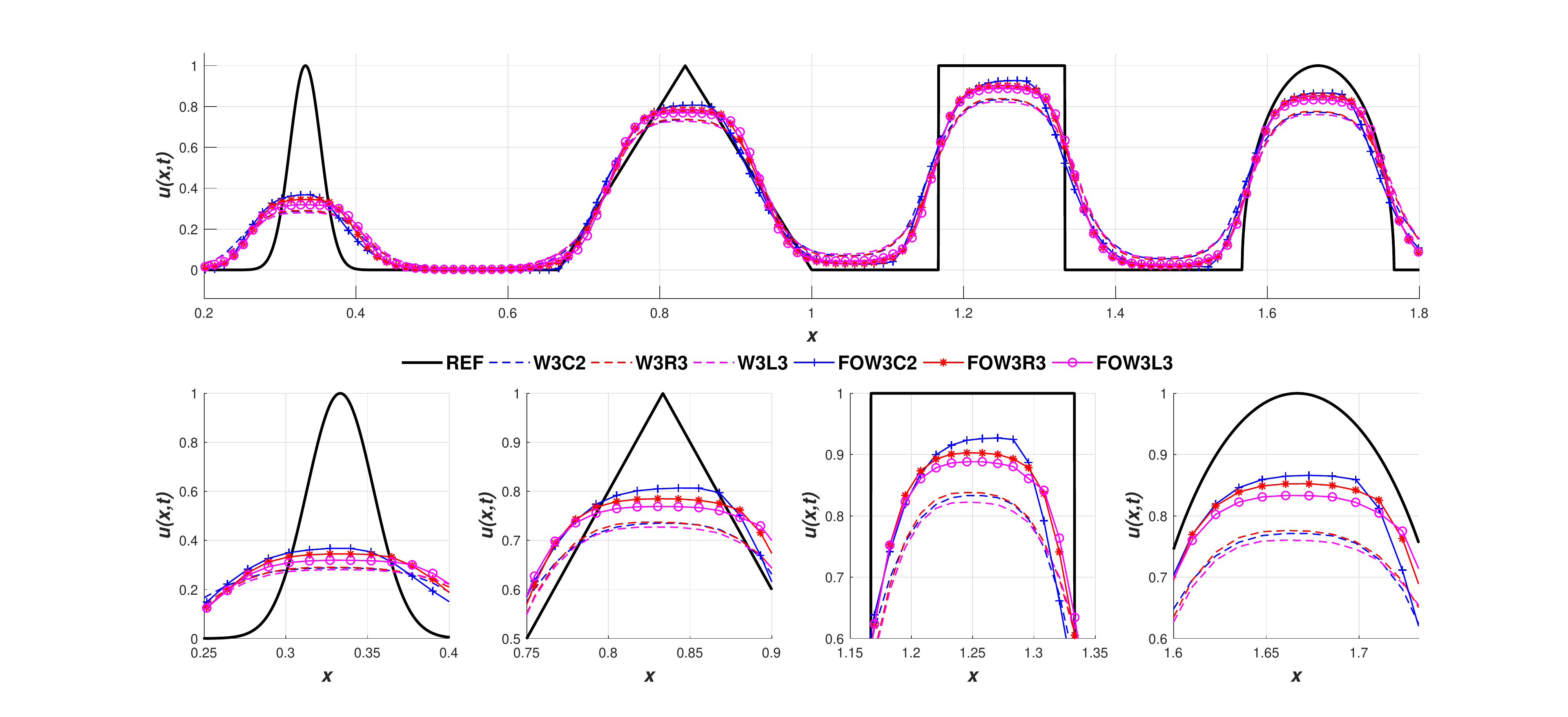}
	\vspace {-1 cm}
	\caption{Test 1: linear transport equation with initial conditions (\ref{initest1}), $CFL=0.5$ and $t=2$s. Methods based on 3rd order reconstructions: general view (up) and zoom of the areas of interest (down).}
	\label{test_1a}
\end{figure}

\begin{figure}[!ht]
   \small
	\setlength{\unitlength}{1mm}
	\centering	
	\includegraphics[width=\textwidth, height=9cm]{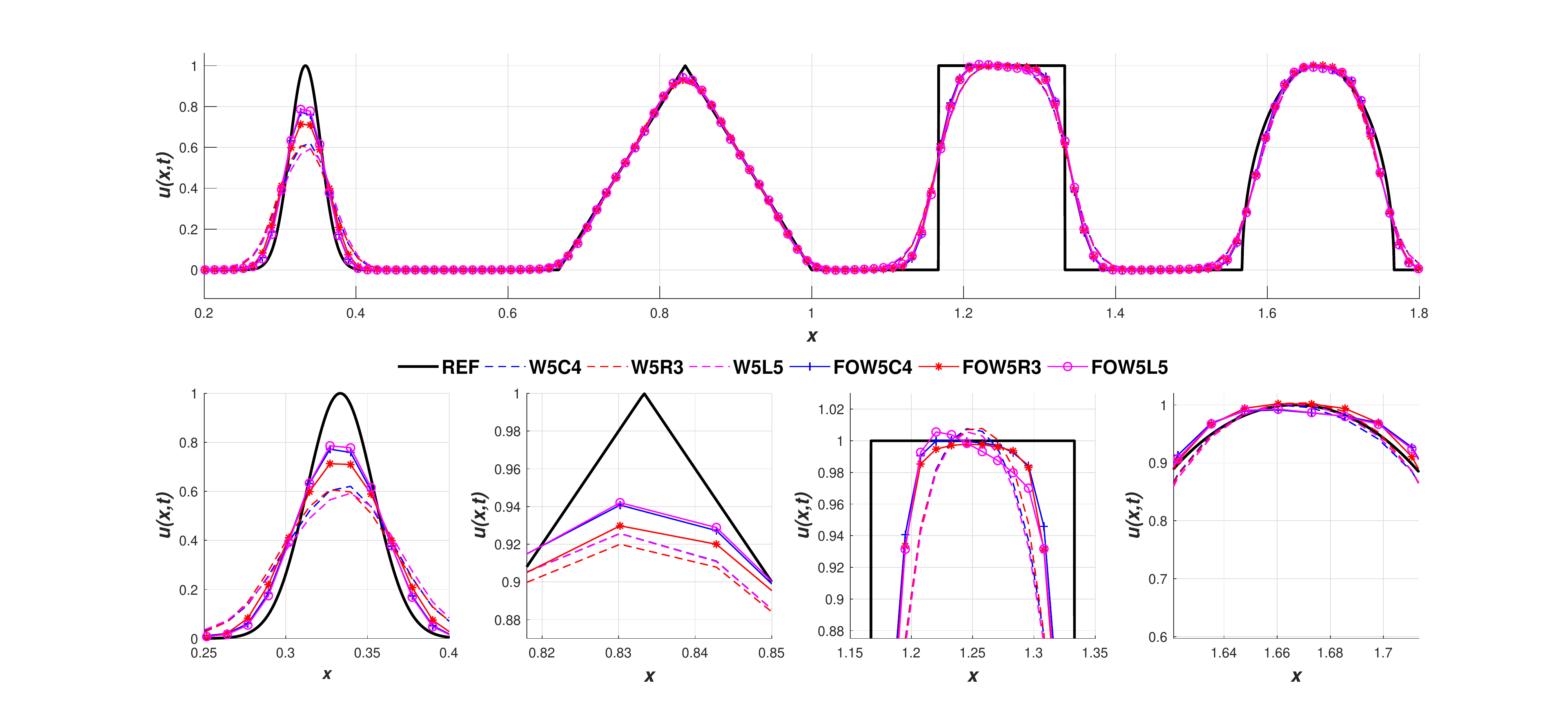}
	\vspace {-1 cm}
	\caption{Test 1: linear transport equation with initial conditions (\ref{initest1}), $CFL=0.5$ and $t=2$s. Methods based on 5th order reconstructions: general view (up) and zoom of the areas of interest (down).}
	\label{test_1b}
\end{figure}

\begin{figure}[!ht]
   \small
	\setlength{\unitlength}{1mm}
	\centering	
	\includegraphics[width=\textwidth, height=9cm]{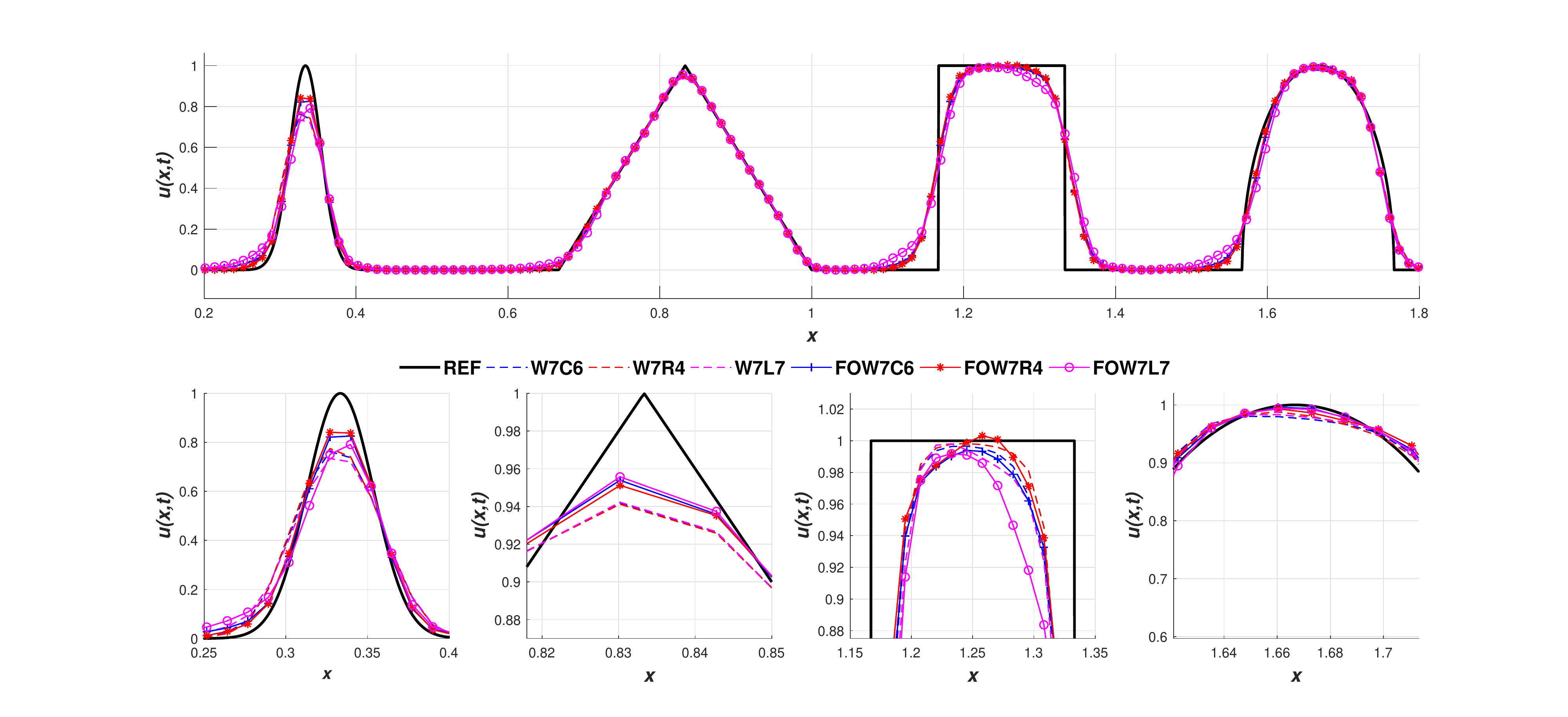}
	\vspace{-1 cm}
	\caption{Test 1: linear transport equation with initial conditions (\ref{initest1}), and $t=2$s. Methods based on 7th order reconstructions with $CFL = 0.5$: general view (up) and zoom of the areas of interest (down).}
	\label{test_1c}
\end{figure}

\begin{figure}[!ht]
   \small
	\setlength{\unitlength}{1mm}
	\centering	
	\includegraphics[width=\textwidth, height=9cm]{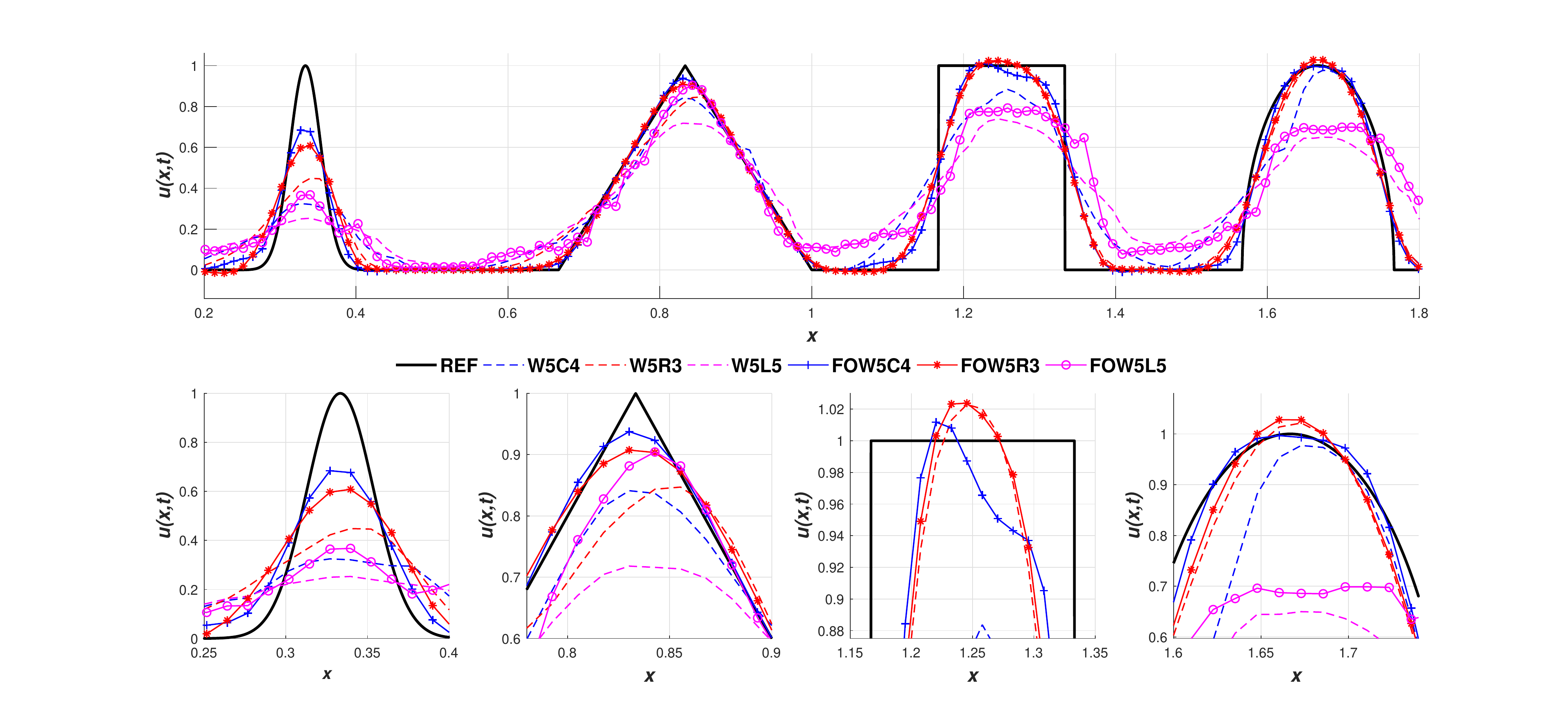}
	\vspace{-1 cm}
	\caption{Test 1: linear transport equation with initial conditions (\ref{initest1}), and $t=2$s. Methods based on 5th order reconstructions with $CFL = 0.9$: general view (up) and zoom of the areas of interest (down).}
	\label{test_1d}
\end{figure}

 Figures \ref{test_1a}, \ref{test_1b}, \ref{test_1c} and  \ref{test_1d} show the results obtained with the methods W3R3, W3C2, W3L3, W5R3, W5C4, W5L5, W7R4, W7C6, W7L7, FOW3R3, FOW3C2, FOW3L3, FOW5R3, FOW5C4, FOW5L5, FOW7R4, FOW7C6,  and FOW7L7 at time $t = 2.$ using a $200$-point mesh, $a=1$, periodic boundary conditions, and $CFL=\{0.5, 0.9\}$. This test is a slight modification of the one proposed by Jiang and Shu in \cite{JiangShu96}. 
 
  From these plots we can conclude:

\noindent For $CFL=0.5$ 

\begin{itemize}
\item Third order reconstructions (Figure \ref{test_1a}): FOWENO reconstructions give better results than WENO reconstructions in all cases. 
We stress the fact that, in spite of its lower order of accuracy, CAT2 gives very good results particularly when combined with FOW3 reconstruction: see enlarged views. 

\item Fifth order reconstructions (Figure \ref{test_1b}):  SSPRK3 gives worse results than CAT4 and LAT5  in the two first areas of interest with both WENO5 and FOWENO5. While CAT4 and LAT5  give similar results when combined with  W5, LAT5 gives better results for FOWENO5: see enlarged views. 

\item Seventh order reconstructions (Figure \ref{test_1c}) : WENO and FOWENO SSPRK4 give solutions that are slightly  better than those given by CAT6 and LAT7.   
\end{itemize}

\noindent For $CFL=0.9$.  
\begin{itemize}

\item Fifth order reconstructions (Figure \ref{test_1d}): LAT5 methods are not stable for this $CFL$ value, and SSPRK4 methods give oscillatory solution, especially near discontinuities. CAT4 combined with  FOWENO5 is stable and gives very good solutions: see enlarged views.
\end{itemize}

Table \ref{times1} shows the CPU times corresponding to the different methods for $t = 2.$ and $CFL = 0.5$. The values are obtained by averaging the computational cost of ten runs. The entries of the table show the ratio between the computational time of each method and the corresponding to W5R3 which is the reference.    

\begin{table}[htbp]
\begin{center}
\resizebox{10 cm}{!}{
\begin{tabular}{|c||c||c||c||c||c|}
\hline
FOW3C2      &   FOW3L3  & FOW3R3  & W3C2   &  W3L3  &  W3R3 \\
 0.3695     &   0.4509  & 0.8351  & 0.3742 & 0.734  &  0.6468 \\
\hline \hline 
FOW5C4      &   FOW5L5  & FOW5R3  &  W5C4 &  W5L5   &  W5R3 \\
1.0546      &  0.7540   & 0.9980  & 1.1936 & 0.7589 &  1  \\
\hline \hline 
FOW7C6      &   FOW7L7  & FOW7R4  & W7C6   &  W7L7  & W7R4 \\
2.5049      &   1.1818  & 4.4116  & 3.4330 & 1.715  & 5.1513 \\
\hline 
\end{tabular}
}
\vspace{2mm}
\caption{CPU time ratios for Test 1: linear transport equation with initial conditions (\ref{initest1}), $CFL = 0.5$, and $t=2$.}
\label{times1}
\end{center}
\end{table}

The following conclusions can be drawn:
\begin{itemize}
    \item The cheapest method is FOW3C2 (that is only second order accurate in time) and the most expensive is W7R4 (due to the extra cost of the smoothness indicators and to the 10 stages of SSPRK4).
    \item Methods based on WENO reconstructions are more costly than their corresponding FOWENO counterparts with the only exception of FOW3R3 and W3R3. Moreover the differences increase with the order.
    \item Methods based on CAT$s$ are more costly than their LAT$(s+1)$ counterparts with the only exception of CAT2. The differences increase with the order. Nevertheless, this extra cost is compensated by the better stability properties of CAT methods for CFL values bigger than 0.5.
\end{itemize}

\subsubsection{Test 2: Burgers equation}
\begin{figure}[!ht]
	\setlength{\unitlength}{1mm}
	\centering	
    \includegraphics[width=\textwidth, height=15 cm]{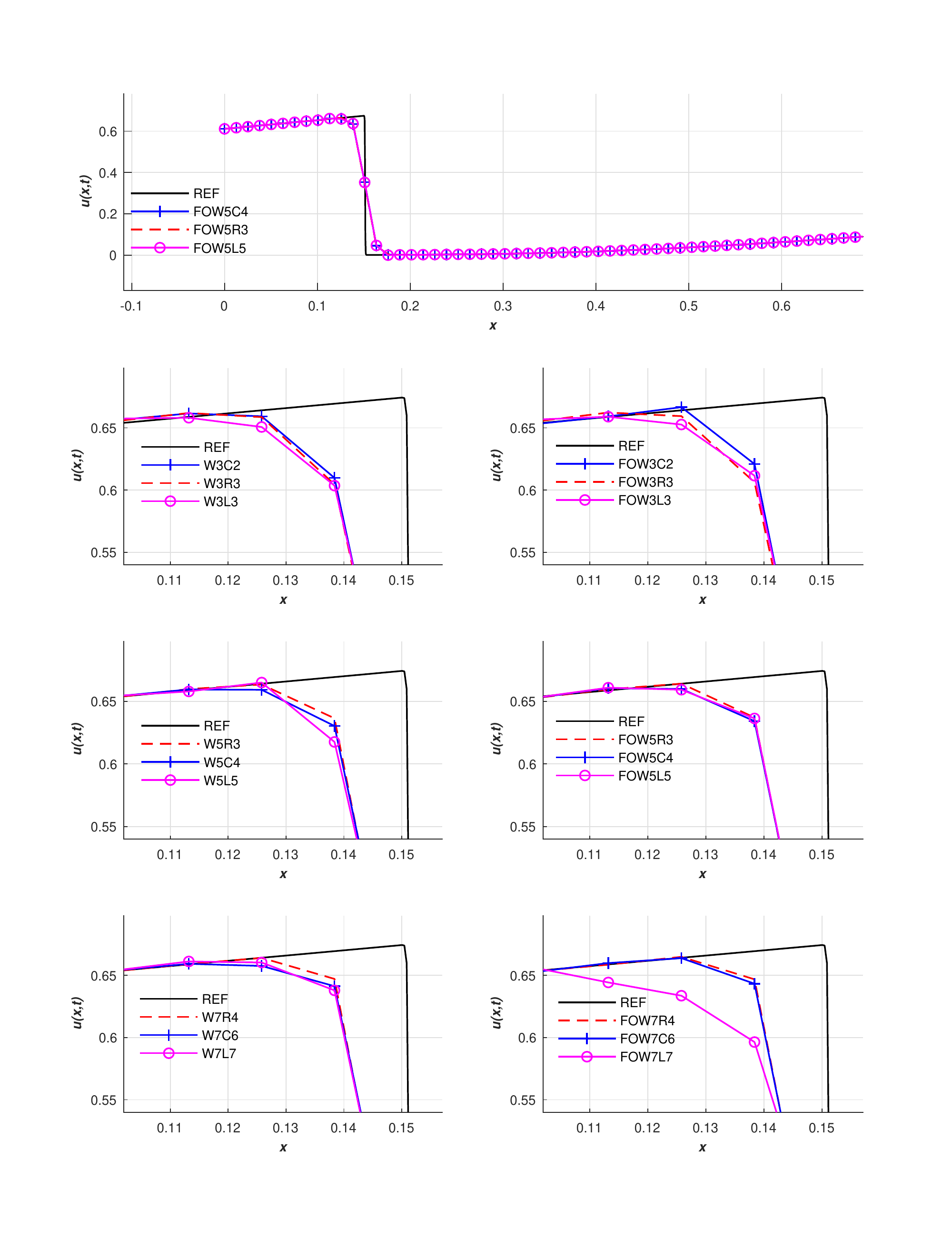}
	\vspace{-1.5 cm}
	\caption{Test 2: Burgers equation with initial conditions (\ref{gauss_distrubution}), $CFL=0.5$ and $t=2$s. 
	Row 1: methods based on 5th order reconstructions: general view. Rows 2-4: zooms of an area of interest.  
	 }
	\label{test2_1}
	\end{figure}

Let us consider now Burgers equation i.e. (\ref{cons_law}) with $f(u)=u^2/2$, in the spatial interval $[0,1]$ with initial condition
\begin{equation} \label{gauss_distrubution} 
u(x,0)= \e^{-10(x-1/2)^2}.
\end{equation}
Figure \ref{test2_1} shows  the numerical solutions obtained with W3R3, W3C2, W3L3, W5R3, W5C4, W5L5, W7R4, W7C6, W7L7, FOW3R3, FOW3C2, FOW3L3,  FOW5R3, FOW5C4, FOW5L5, FOW7R4, FOW7C6 and FOW7L7 methods using a $160$-point mesh, periodic boundary conditions, $CFL=0.5$, and $t=2$s. The numerical results are shown in groups of three to facilitate the comparisons.   From the enlarged views (close to  the shock) the following conclusions can be drawn:

\begin{itemize}
\item Methods based on third order reconstructions (Figure \ref{test2_1} row 2): all the methods based on WENO3 give essentially the same solutions. Some improvements are achieved with FOWENO3 and CAT2 is slightly sharper than the rest.

\item Methods based on fifth order reconstructions (Figure \ref{test2_1} row 3): the results are better than  the ones corresponding to third-order reconstructions as expected. There are no big differences between them, but a slight improvement can be observed when FOWENO reconstructions are used.

\item Methods based on seventh order reconstructions (Figure \ref{test2_1} row 4): WENO7 and FOWENO7 reconstructions give non-oscillatory solutions and better results than third or fifth order reconstructions for CAT6 and RK4, which is not the case for LAT7. 

\end{itemize}

Concerning the quality of the numerical results with $CFL = 0.9$ or the computational cost, the conclusions are similar to the previous test case.

\subsection{1D Systems of conservation laws} 
We consider the 1D Euler equations of gas dynamics:
\begin{equation}
\mathbf{w}_{t}+\mathbf{f}(\mathbf{w})_{x}=\mathbf{0} \text { , }
\end{equation}
where
$$
\mathbf{w}=\left(\begin{array}{c}{\rho} \\ {\rho u} \\ {E}\end{array}\right), \quad \mathbf{f}(\mathbf{w})=\left(\begin{array}{c}{\rho u} \\ {\rho u^2 + p} \\ {u(E+p)}\end{array}\right).
$$
\\
Here, $\rho$ is the density, $u$ the velocity, $E$  the total energy per unit volume and $p$  the pressure. We assume an ideal gas with the equation of state 
\begin{equation*}
p(\rho, e) = (\gamma - 1)\rho e,
\end{equation*}
where $\gamma$ is the ratio of specific heat capacities of the gas and $e$ the internal energy per unit mass given by:
\begin{equation*}
E(\rho, u, e) = \rho (e+\frac{1}{2} u^2).
\end{equation*}

We consider the following 1D Riemann problems whose data are given in Table \ref{torotable}:
\begin{itemize}
\item {\textbf{Test 3}}: Sod  problem \cite{Sod1978}. The solution consists of a left rarefaction, a left contact and a right shock. 

\item {\textbf{Test 4}}: 123 Einfeldt \cite{Einfeldt1991}. The solution consists of two strong rarefactions and a stationary contact discontinuity. The pressure $p$ is  small (close to vacuum).

\item {\textbf{Test 5}}: left half of the blast wave problem \cite{Wood1984}. The solution contains a left rarefaction, a contact and a right shock. 

\item {\textbf{Test 6}}: right half of the blast wave problem \cite{Wood1984}. The solution contains a left shock, a contact discontinuity and a right rarefaction.

\item {\textbf{Test 7}}: blast wave problem \cite{Wood1984}. The solution represents the collision of the right and left shocks corresponding to tests 3 and 4, and consists of a left facing shock (travelling very slowly to the right), a right contact discontinuity and a right shock wave. 
\end{itemize}

The equations are solved in the spatial domain $x \in [0,1]$ with outflow-inflow boundary conditions and a $200$-point mesh. $CFL = 0.9, 0.5, 0.25$ are used for methods based on  with 3rd, 5th, and 7th order reconstructions respectively. We consider WENO reconstructions with $\eps=1e-6$ as in \cite{Shu1989} and FOWENO reconstructions with $\eps=1e-100$ as in \cite{ZBBM2019}. The numerical solutions are compared against the exact solution provided by the HE-E1RPEXACT solver introduced in \cite{Toro2009book}

\begin{table}[htbp]
\begin{center}
\resizebox{\textwidth}{!}{
\begin{tabular}{|c|c|c|c|c|c|c|c|}
\hline
Test  &$\rho_{\mathrm{L}}$  & $u_{\mathrm{L}}$ & $p_{\mathrm{L}}$ & $\rho_{\mathrm{R}}$ & $u_{\mathrm{R}}$ & $p_{\mathrm{R}}$ & time  (sec.) \\
\hline 
\textbf{3} & 1.0 & 0.0 & 1.0 & 0.125 & 0.0 & 0.1 & {0.25} \\ 
\hline
\textbf{4} & {1.0} & {-2.0} & {0.4} & {1.0} & {2.0} & {0.4} & {0.15} \\
\hline
 \textbf{5} & {1.0} & {0.0} & {1000.0} & {1.0} & {0.0} & {0.01} & {0.012}\\ 
\hline
 \textbf{6} & {1.0} & {0.0} & {0.01} & {1.0} & {0.0} & {100.0} & {0.035}\\
\hline
 \textbf{7} & {.99924} & {19.5975} & {460.894} & {5.99242} & {-6.19633} & {46.0950} & {0.035} \\
\hline
\end{tabular}
}
\vspace{2mm}
\caption{Riemann problems for 1D Euler equations.}
\label{torotable}
\end{center}
\end{table}

\begin{figure}[!ht]
	\small
	\setlength{\unitlength}{1mm}
	\centering	
	\includegraphics[width=\textwidth,height=18 cm]{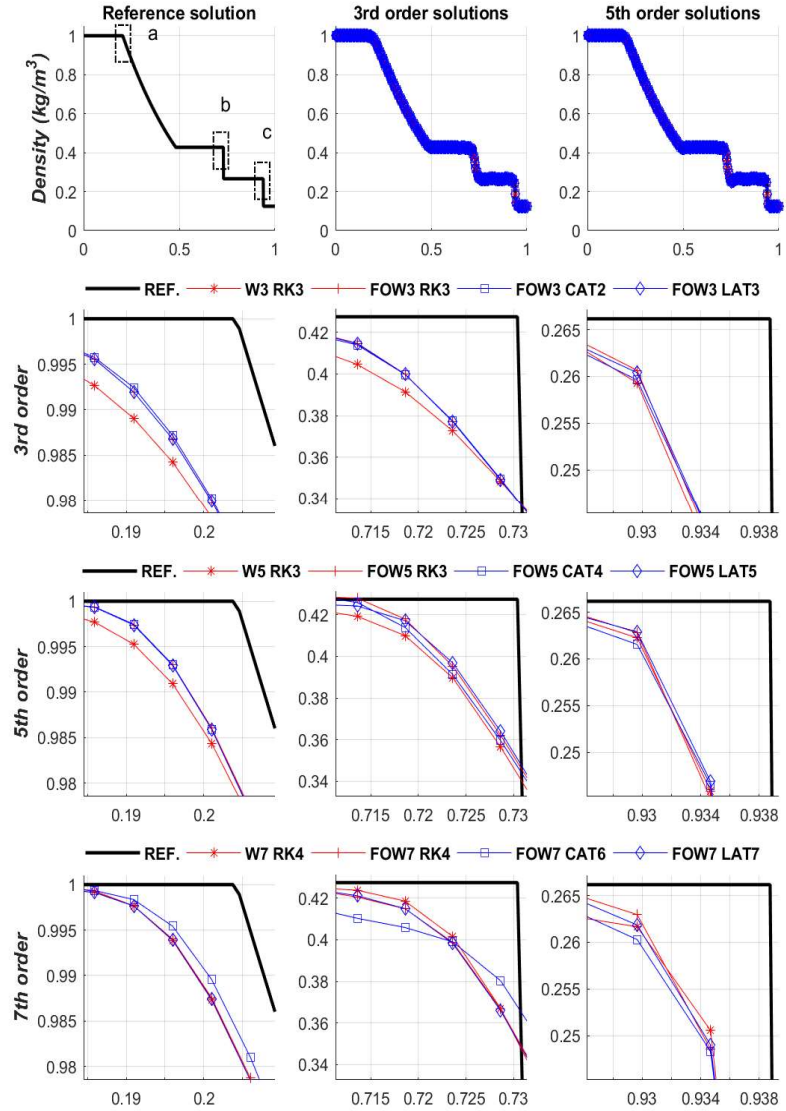}
	\vspace{-2 cm}
	\caption{ Test 3:  1D Euler equations. Sod problem: density. Row 1: exact solution (left), methods using 3rd order (center)  and  5th order (right) reconstruction operators. Rows 2-4: zooms corresponding to areas  \textit{a}, \textit{b} and \textit{c}. $CFL = 0.9, 0.5, 0.25$  for methods based on  with 3rd, 5th, and 7th order  reconstructions respectively.}
	\label{1d_sod_1}
	\end{figure}

\begin{figure}[!ht]
	\small
	\setlength{\unitlength}{1mm}
	\centering	
	\includegraphics[width=\textwidth,height=15 cm]{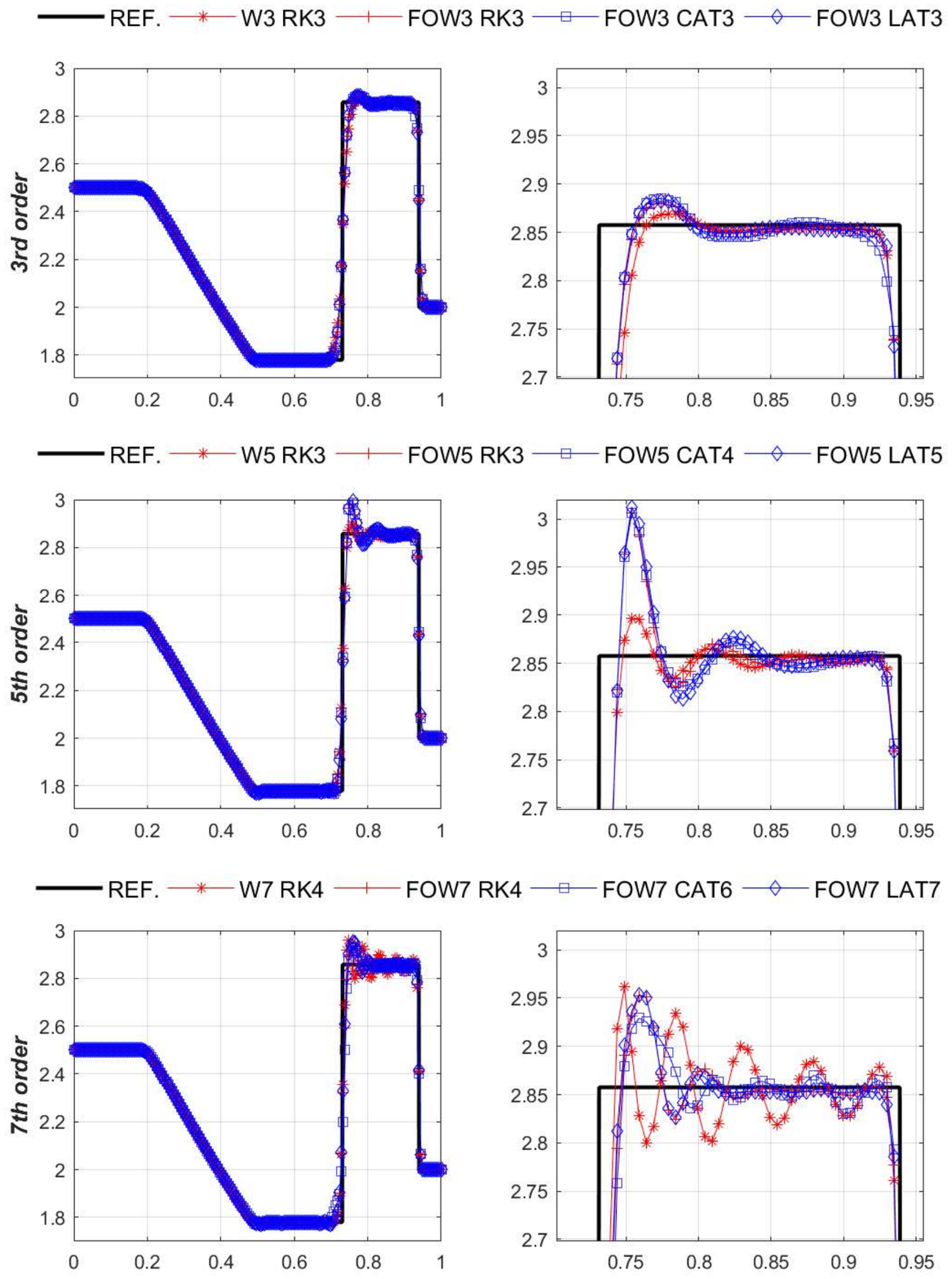}
	\vspace{-1.7 cm}
	\caption{ Test 3:  1D Euler equations. Sod problem: internal energy. Methods using 3rd order (row 1), 5th order (row 2), and 7th order (row 3) reconstruction operators. Left: general  view. Right:  zoom of an area of interest.  Exact solution: black line.  $CFL = 0.9, 0.5, 0.25$  for methods based on  with 3rd, 5th, and 7th order  reconstructions respectively.}
	\label{1d_sod_2}
	\end{figure}

\begin{figure}[!ht]
	\small
	\setlength{\unitlength}{1mm}
	\centering	
	\includegraphics[width=\textwidth,height=16 cm]{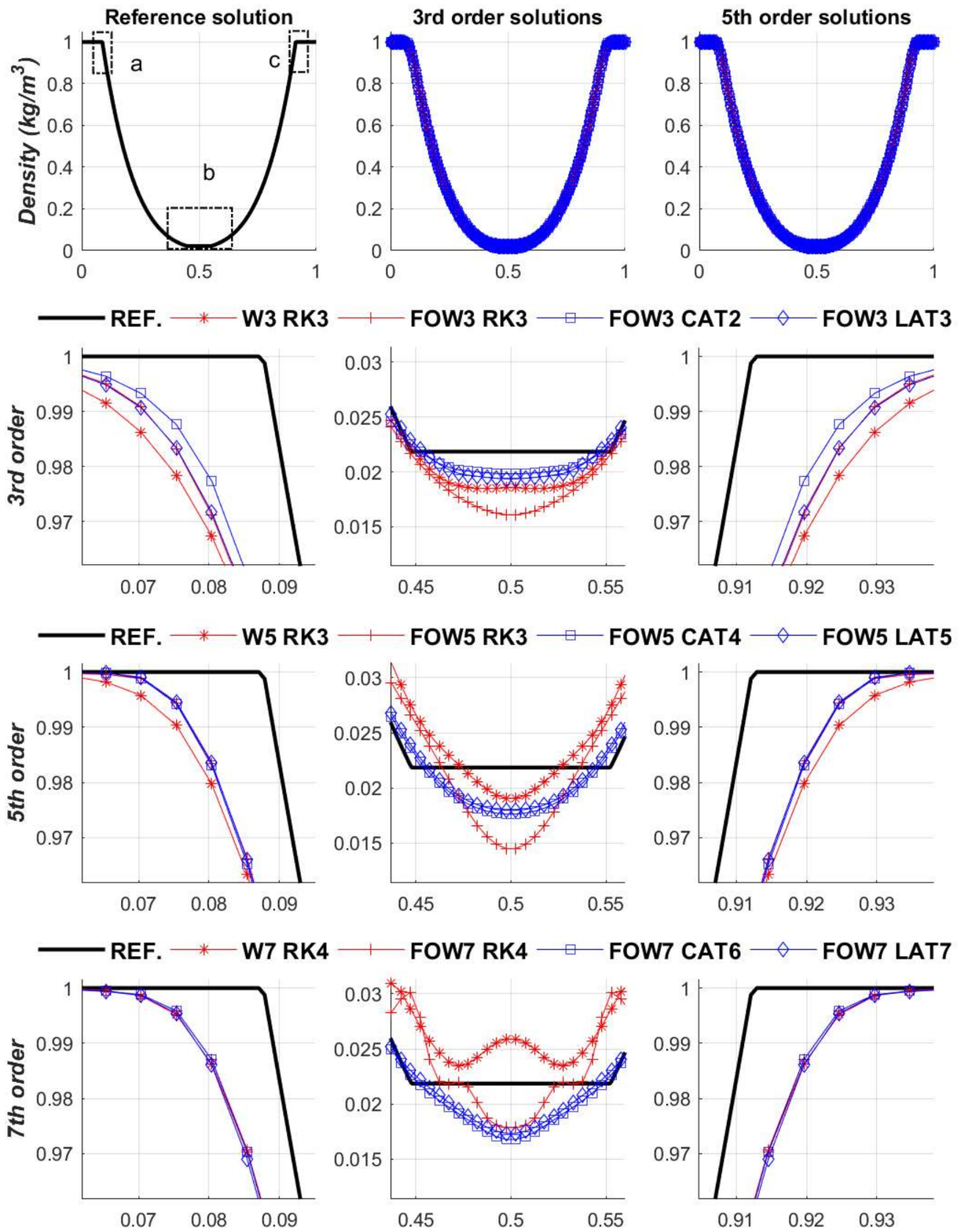}
	\vspace{-1.7cm}
	\caption{Test 4:  1D Euler equations. 123 Einfeldt problem: density.  Row 1: exact solution (left), methods using 3rd (center)  and  5th order (right) reconstruction operators. Rows 2-4: zooms corresponding to areas  \textit{a}, \textit{b} and \textit{c}. $CFL = 0.9, 0.5, 0.25$  for methods based on  with 3rd, 5th, and 7th order  reconstructions respectively.}
	\label{1d_123_1}
	\end{figure}

\begin{figure}[!ht]
	\small
	\setlength{\unitlength}{1mm}
	\centering	
	\includegraphics[width=\textwidth, height=14 cm]{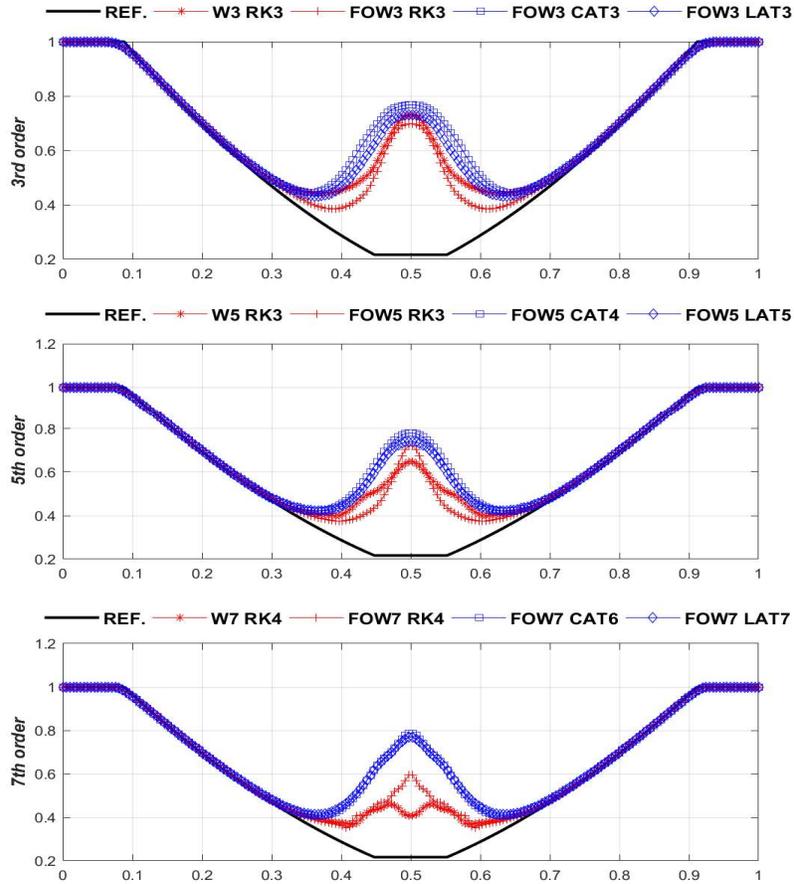}
	\vspace{-1.5 cm}
	\caption{Test 4:  1D Euler equations. 123 Einfeldt problem: internal energy. Methods using 3rd order (row 1), 5th order (row 2), and 7th order (row 3) reconstruction operators.  Exact solution: black line.  $CFL = 0.9, 0.5, 0.25$  for methods based on  with 3rd, 5th, and 7th order  reconstructions respectively.}
	\label{1d_123_2}
	\end{figure}

\begin{figure}[!ht]
	\small
	\setlength{\unitlength}{1mm}
	\centering	
	\includegraphics[width=\textwidth, height=18 cm]{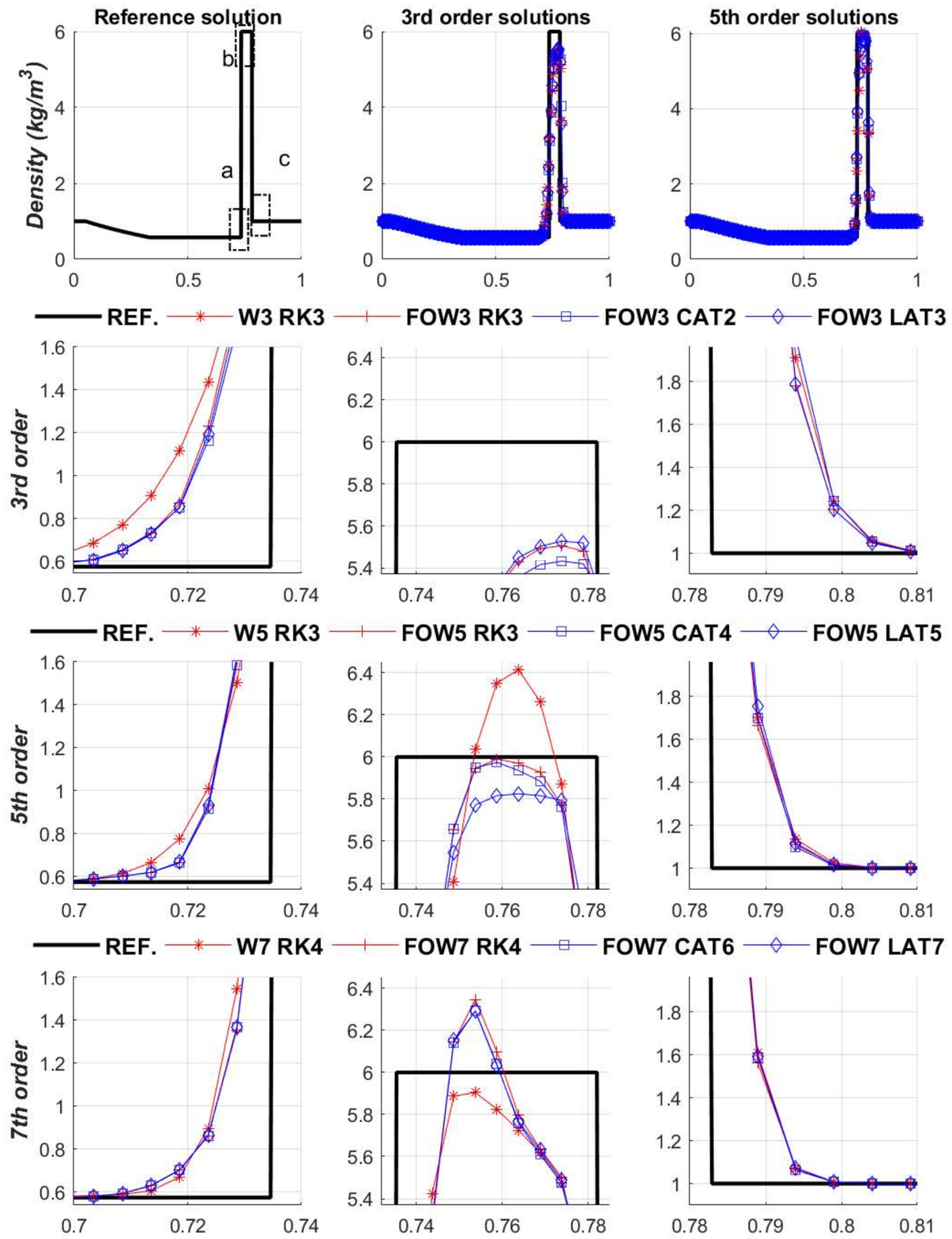}
	\vspace{-1.7 cm}
	\caption{  Test 5:  1D Euler equations. Left half of the blast wave problem of Woodward and Colella: density. Row 1: exact solution (left), methods using 3rd (center)  and  5th order (right) reconstruction operators. Rows 2-4: zooms corresponding to areas  \textit{a}, \textit{b} and \textit{c}.  $CFL = 0.9, 0.5, 0.25$  for methods based on  with 3rd, 5th, and 7th order  reconstructions respectively.}
	\label{1d_wc_r_1}
	\end{figure}

\begin{figure}[!ht]
	\small
	\setlength{\unitlength}{1mm}
	\centering	
	\includegraphics[width=\textwidth, height=14 cm]{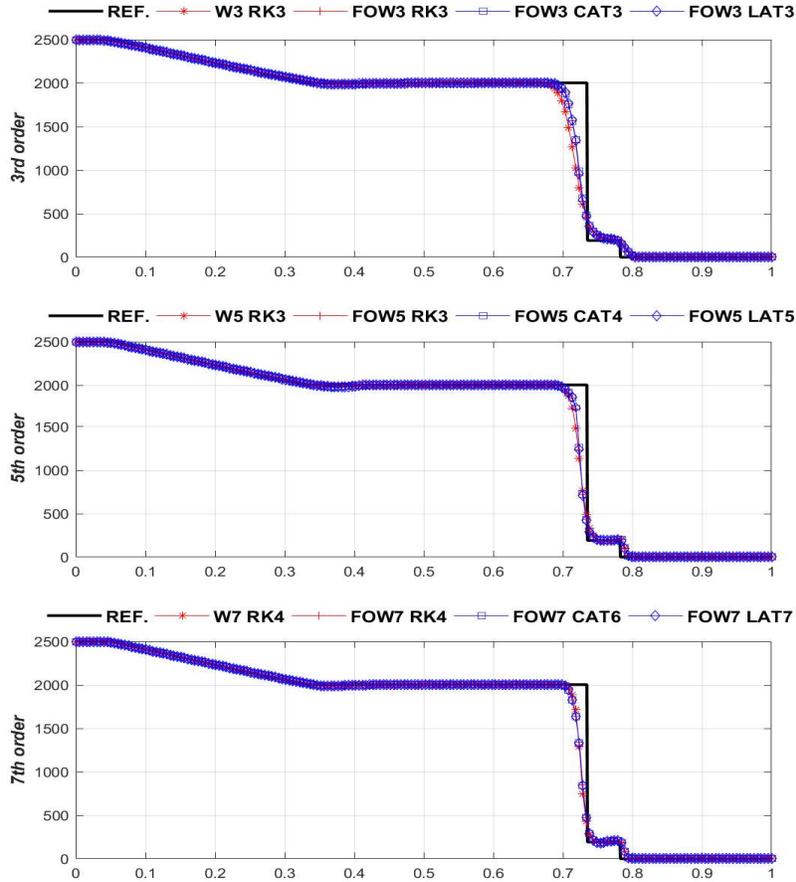}
	\vspace{-1.7 cm}
	\caption{Test 5:  1D Euler equations. Left half of the blast wave problem of Woodward and Colella: internal energy. 
	Methods using 3d order (row 1), 5th order (row 2), and 7th order (row 3) reconstruction operators. Exact solution: black line.  $CFL = 0.9, 0.5, 0.25$  for methods based on  with 3rd, 5th, and 7th order  reconstructions respectively.}
	\label{1d_wc_r_2}
	\end{figure}

\begin{figure}[!ht]
	\small
	\setlength{\unitlength}{1mm}
	\centering	
	\includegraphics[width=\textwidth, height=16 cm]{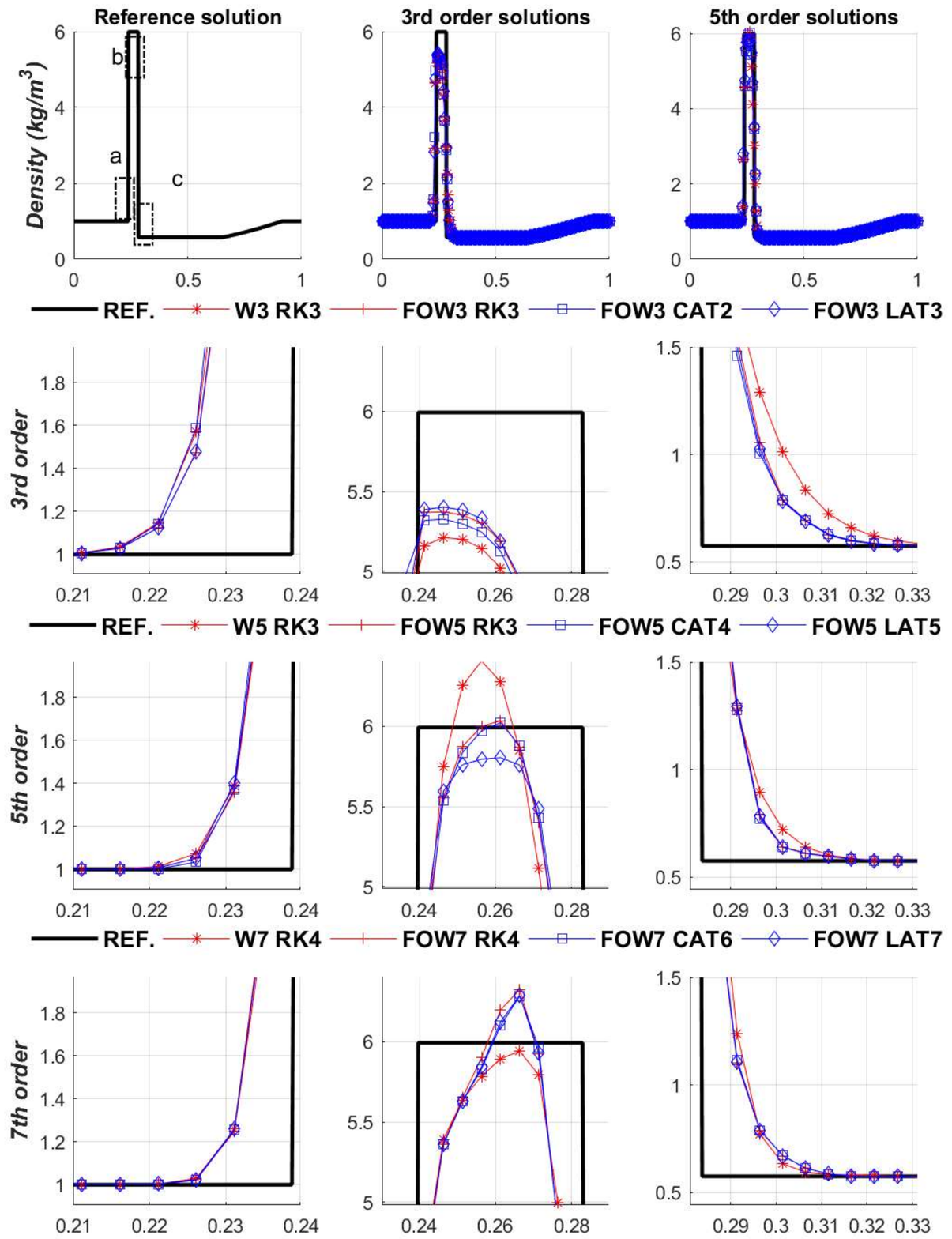}
	\vspace{-1.5 cm}
	\caption{ Test 6:  1D Euler equations. Right half of the blast wave problem of Woodward and Colella: density. Row 1: exact solution (left), methods using 3rd order (center)  and  5th order (right) reconstruction operators. Rows 2-4: zooms corresponding to areas  \textit{a}, \textit{b} and \textit{c}.  $CFL = 0.9, 0.5, 0.25$  for methods based on  with 3rd, 5th, and 7th order  reconstructions respectively.}
	\label{1d_wc_l_1}
	\end{figure}

\begin{figure}[!ht]
	\small
	\setlength{\unitlength}{1mm}
	\centering	
\includegraphics[width=\textwidth, height=14 cm]{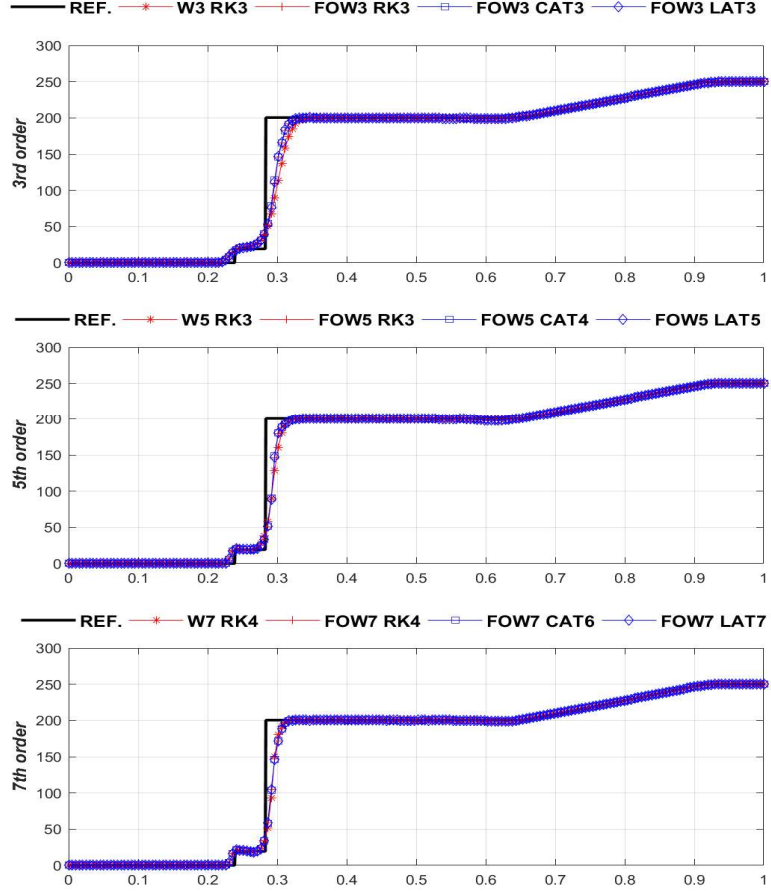}
	\vspace{-1.5 cm}
	\caption{Test 6:  1D Euler equations. Right half of the blast wave problem of Woodward and Colella: internal energy.
	Methods using 3rd order (row 1), 5th order (row 2) and 7th order (row 3) reconstruction operators. Left: general view. Right: zoom of an area of interest.  Exact solution: black line.  $CFL = 0.9, 0.5, 0.25$  for methods based on  with 3rd, 5th, and 7th order  reconstructions respectively.}
	\label{1d_wc_l_2}
	\end{figure}

\begin{figure}[!ht]
	\small
	\setlength{\unitlength}{1mm}
	\centering	
	\includegraphics[width=\textwidth, height=16 cm]{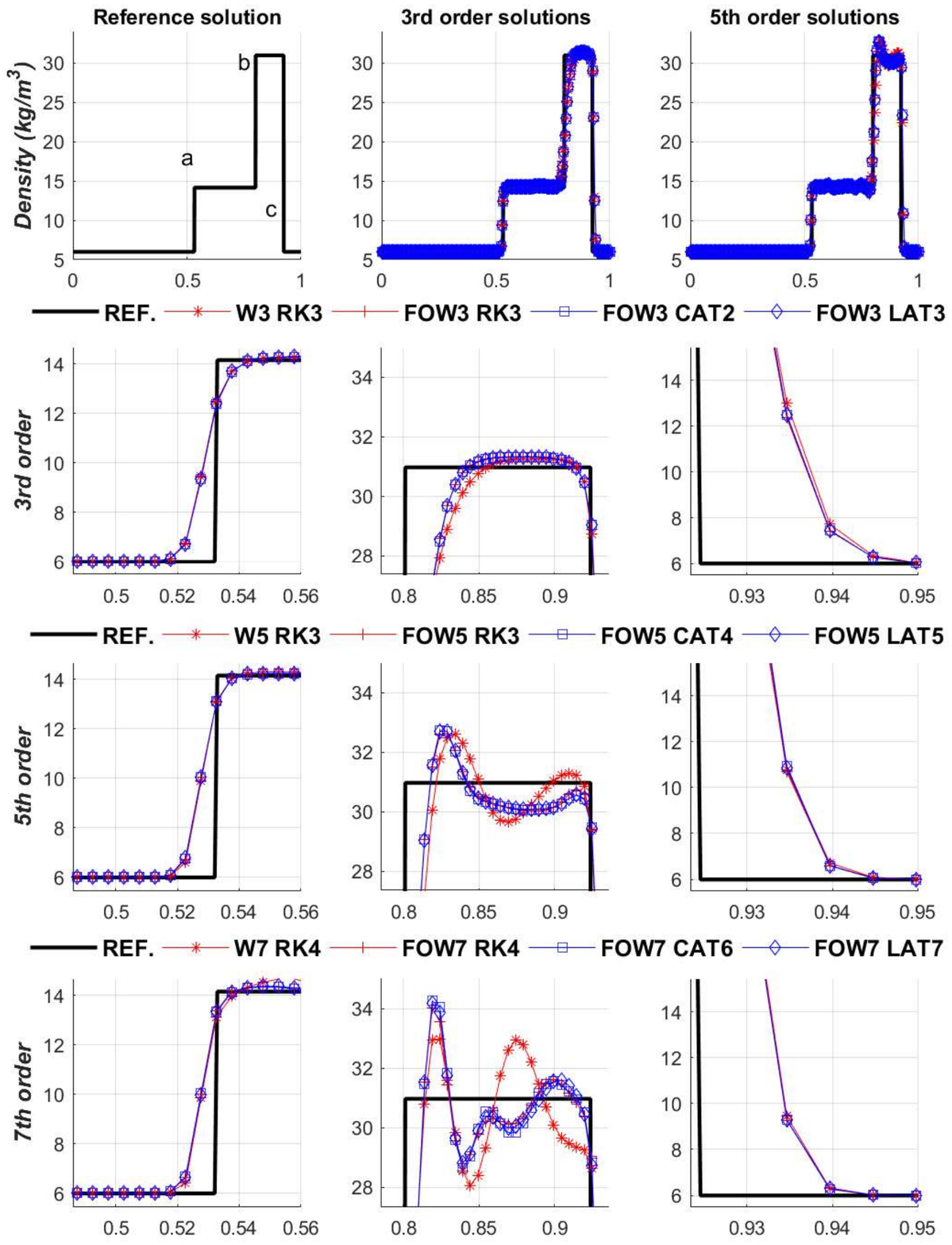}
	\vspace{-1.5 cm}
	\caption{ Test 7:  1D Euler equations. Woodward and Colella problem: density. Row 1: exact solution (left),  methods using 3rd order (center)  and  5th order (right) reconstruction operators. Rows 2-4: zooms corresponding to areas  \textit{a}, \textit{b} and \textit{c}.  $CFL = 0.9, 0.5, 0.25$  for methods based on  with 3rd, 5th, and 7th order  reconstructions respectively.}
	\label{1d_wc_1}
	\end{figure}

\begin{figure}[!ht]
	\small
	\setlength{\unitlength}{1mm}
	\centering	
	\includegraphics[width=\textwidth, height=14 cm]{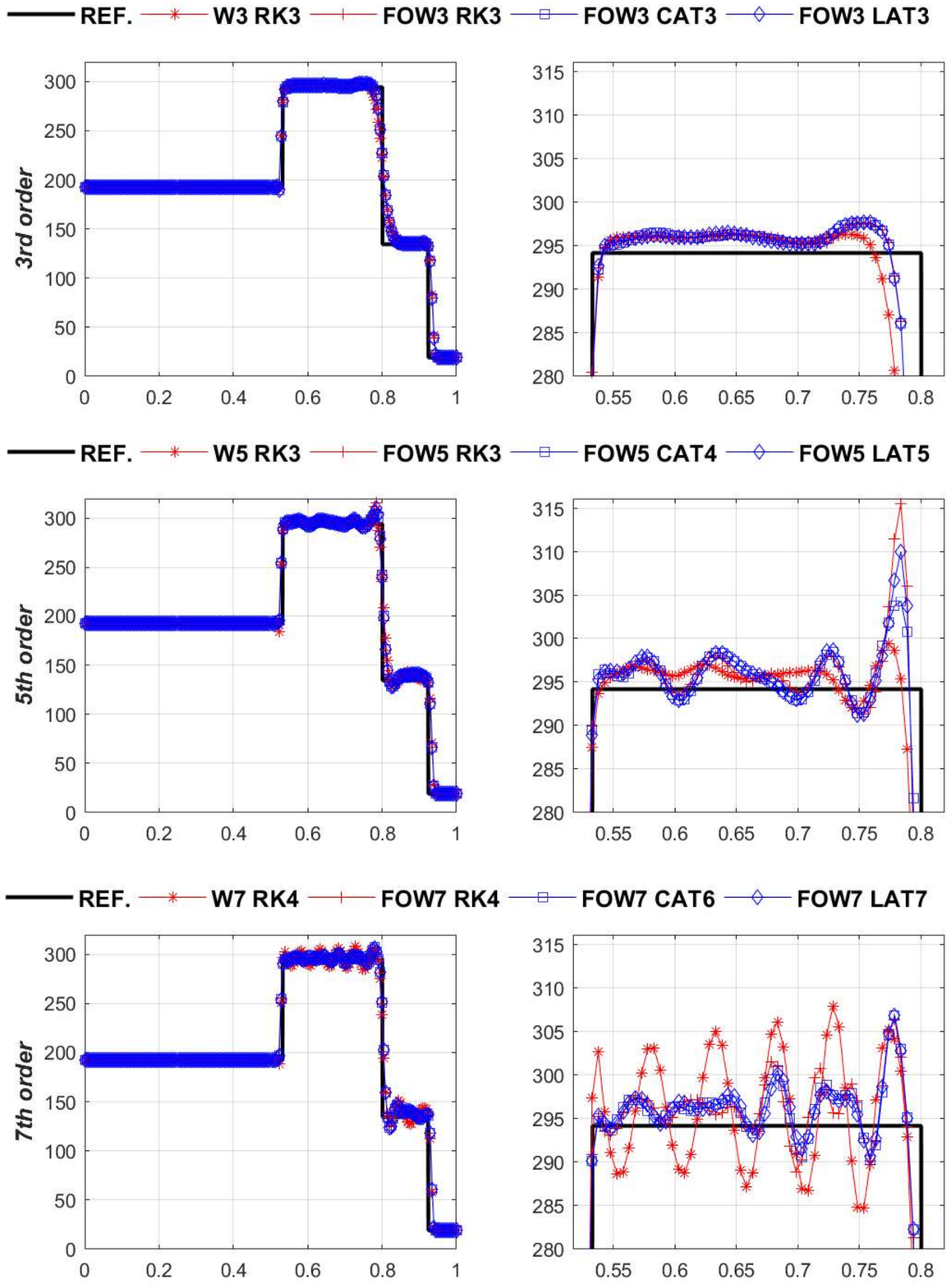}
	\vspace{-1.5 cm}
	\caption{Test 7:  1D Euler equations. Woodward and Colella problem: internal energy. Methods using 3rd order (row 1), 5th order (row 2) and 7th order (row 3) reconstruction operators. Left: general view. Right: zoom of an area of interest. Exact solution: black line.  $CFL = 0.9, 0.5, 0.25$  for methods based on  with 3rd, 5th, and 7th order  reconstructions respectively.}
	\label{1d_wc_2}
	\end{figure}

The numerical results are shown in Figures \ref{1d_sod_1}-\ref{1d_wc_2}.
Two figures are shown for every test case, the first one corresponds to densities and the second one to internal energies. In the first row of the figures corresponding to the densities, we show the global views of the reference and the numerical solutions obtained using third and fifth order reconstructions. Rows 2-4 show enlarged views of the areas of interest labelled \textit{a}, \textit{b} and \textit{c} in the global view of the reference solution. In the figures corresponding to internal energy we plot global views of the numerical results for third, fifth and seventh order reconstructions (left column) and enlarged views of an interest area of each one of them (right column).

\begin{table}[htbp]
\begin{center}
\resizebox{10 cm}{!}{
\begin{tabular}{|c||c||c||c||c||c|}
\hline
FOW3C2      &   FOW3L3  & FOW3R3  & W3C2   &  W3L3  &  W3R3 \\
1.1830      &   1.6352   & 2.8026  & 1.0000 & 1.3744  &  2.1764 \\
\hline \hline 
FOW5C4      &   FOW5L5  & FOW5R3  &  W5C4    &  W5L5   &  W5R3 \\
5.0546      &   3.4400   & 3.2980  &  5.1642 &  3.7589 &  3.5268  \\
\hline \hline 
FOW7C6      &   FOW7L7  & FOW7R4  & W7C6    &  W7L7      & W7R4 \\
23.8827      &   18.1818  &  19.7516  & 29.9430 &  22.7150  & 29.9490 \\
\hline 
\end{tabular}
}
\vspace{2mm}
\caption{CPU time ratios for test 7: 1D Euler equations with the Woodward and Colella problem, $CFL = 0.25$, and $t=0.035$s.}
\label{times2}
\end{center}
\end{table}

\begin{itemize}

\item{\textbf{Test 3}}: Figures  \ref{1d_sod_1} and \ref{1d_sod_2}. In general, all the solutions are acceptable and their quality improve with the order of accuracy. Methods based on  FOWENO reconstruction are slightly sharper than those based on WENO with exception of FOW7C6 near the contact discontinuity (the approximation obtained of this wave is worse but oscillations appear, even for long-time simulation). Concerning the internal energies, solutions obtained with LAT and CAT are less oscillatory: see the enlarged views.  

\item{\textbf{Test 4}}: Figures \ref{1d_123_1} and \ref{1d_123_2}. This is a hard test in which significant differences between WENO and FOWENO reconstructions can be seen. For densities, FOW3C2 and FOW3L3 give the closest solutions to the reference in area \textit{b}. Moreover, all FOWENO-AT solutions are stable and non-oscillatory. For internal energies,  solutions corresponding to WENO methods show oscillations but they are closer to the exact solution.   

\item{\textbf{Test 5}}:  Figures
  \ref{1d_wc_r_1} and \ref{1d_wc_r_2}. 3rd order accuracy is not enough in this case to capture good solutions, especially in  area \textit{c}.  FOW5CAT4 and FOW5LAT5 give better solutions than W5R3, which is under dissipative. However, for  seventh order reconstruction the situation is the opposite,  due to the use of SSPRK$\_10\_4$ for WENO7. For internal energies, no significant differences are detected.

\item{\textbf{Test 6}}: Figures \ref{1d_wc_l_1} and \ref{1d_wc_l_2}. Similar conclusions to Test 5.

\item{\textbf{Test 7}}:  Figures 
  \ref{1d_wc_1} and \ref{1d_wc_2}. In order to compare the cpu times, $CFL=0.25$ has been chosen for all the methods.  Methods based on 7th order reconstructions  give the best approximations in areas \textit{a} and \textit{c} but produce some oscillations in  area \textit{b}. These oscillations are particularly noticeable in the top part of the internal energy  solutions, in which the solutions provided by AT methods are less oscillatory.  CPU times are shown in Table \ref{times2}.  WENO3-CAT2 (which is the faster method) is the reference. Some conclusions can be drawn from this table:
  
  \begin{enumerate}
      \item 3rd order methods based on WENO are cheaper than FOMENO3: in this case the smooth indicators are the same and FOWENO has the extra computational cost due to the computation of the optimal weights.
    
    \item For reconstructions of order 5th or greater, methods based on  FOWENO are faster than those based on WENO.
      
      \item To pass from C2 to C4 using the same reconstruction operator multiplies the computational time approximately by 3. And to pass from C4 to C6 by a factor between 4 and 6. 
      
      \item To pass from L3 to L5 using the same reconstruction operator multiplies the computational time approximately by 5. And to pass from L5 to L7 by a factor between 6 and 7. 
      
      \item To pass from R2 to R3 using the same reconstruction operator multiplies the computational time approximately by 1.5. And to pass from R3 from R5 by a factor between 6 and 8.5.

  \end{enumerate}
  
  \end{itemize}


\subsection{2D Systems of conservation laws}

We consider now the two-dimensional Euler equations of gas dynamics:
\begin{equation}
\mathbf{w}_{t}+\mathbf{f}(\mathbf{w})_{x} +\mathbf{g}(\mathbf{w})_{y}=\mathbf{0} \text { , }
\end{equation}
where
$$
\mathbf{w}=\left(\begin{array}{c}{\rho} \\ {\rho u} \\ {\rho v} \\ {E}\end{array}\right), 
\quad \mathbf{f}(\mathbf{w})=\left(\begin{array}{c}{\rho u} \\ {\rho u^2 + p} \\ {\rho u v}  \\ {u(E+p)}\end{array}\right) , 
\quad \mathbf{g}(\mathbf{w})=\left(\begin{array}{c}{\rho v} \\ {\rho u v} \\ {\rho v^2  + p} \\ {v(E+p)}\end{array}\right)  .
$$
\\
$\rho$ is again the density; $u,v$ are the components of the velocities in the $x,y$ directions respectively; $E$,  the total energy per unit volume; and $p$, the pressure. The  equation of state 
\begin{equation}
p(\rho,u,v,E) = (\gamma - 1)\left(E-\frac{\rho}{2}(u^2+v^2)\right),
\end{equation}
is assumed again where  $\gamma$ is the ratio of specific heat capacities of the gas.

\noindent From the nineteen configurations of the 2-D Riemann problems presented in \cite{Lax1998} six relevant configurations have been selected, namely:  3, 6, 11, 13, 17 and 19. The initial data of the Riemann problems consist of constant states at every quadrant of the spatial domain that are chosen so that the 1D Riemann problems corresponding to two adjacent states consist of only one one-dimensional simple wave: a shock $\textit{S}$, a rarefaction wave $\textit{R}$, or  a  slip line i.e. a contact discontinuity with discontinuous tangential velocity $\textit{J}$.  The sub-indexes  $(l,r) \in \{ (2,1),(3,2),(3,4),(4,1) \}$  indicate the involved quadrants. For shock and rarefactions an over-arrow indicate the direction  (backward or forward). And for contact discontinuities a sign $+/-$ is used (instead of the over-arrow), to denote whether it is a positive or negative slip line. Full information and analysis  can be found in \cite{Lax1998}.

\begin{figure}[!ht]
	\small
	\setlength{\unitlength}{1mm}
	\centering
	\includegraphics[width=\textwidth, height=19 cm]{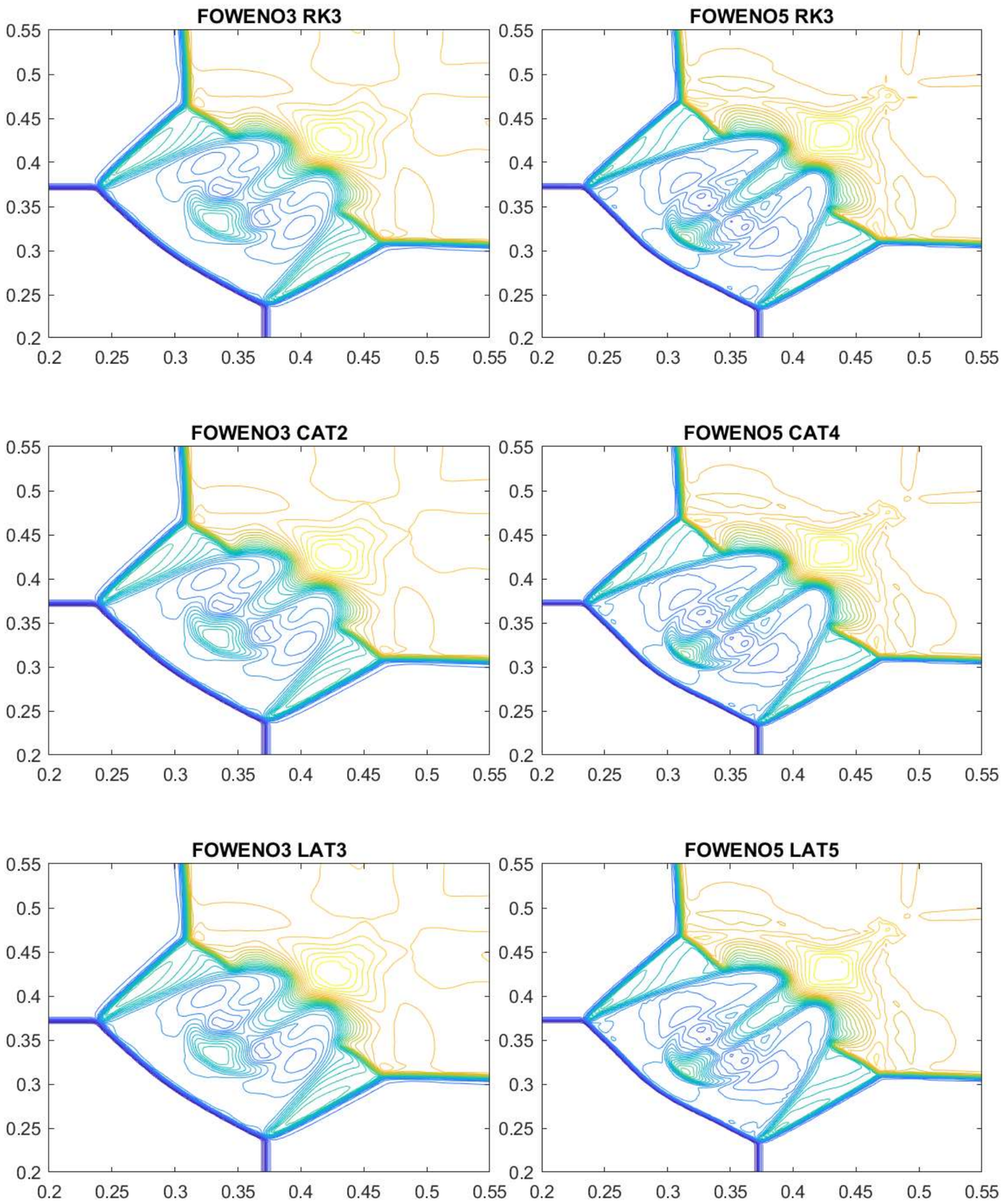}
	\vspace{ -1 cm}
	\caption{Test 8: 2D Euler equations. Lax configuration 3: density computed with  FOWENO-RK, FOWENO-CAT and  FOWENO-LAT.}
	\label{2d_test_3}
	\end{figure}

\begin{figure}[!ht]
	\small
	\setlength{\unitlength}{1mm}
	\centering	
	\includegraphics[width=\textwidth, height=19 cm]{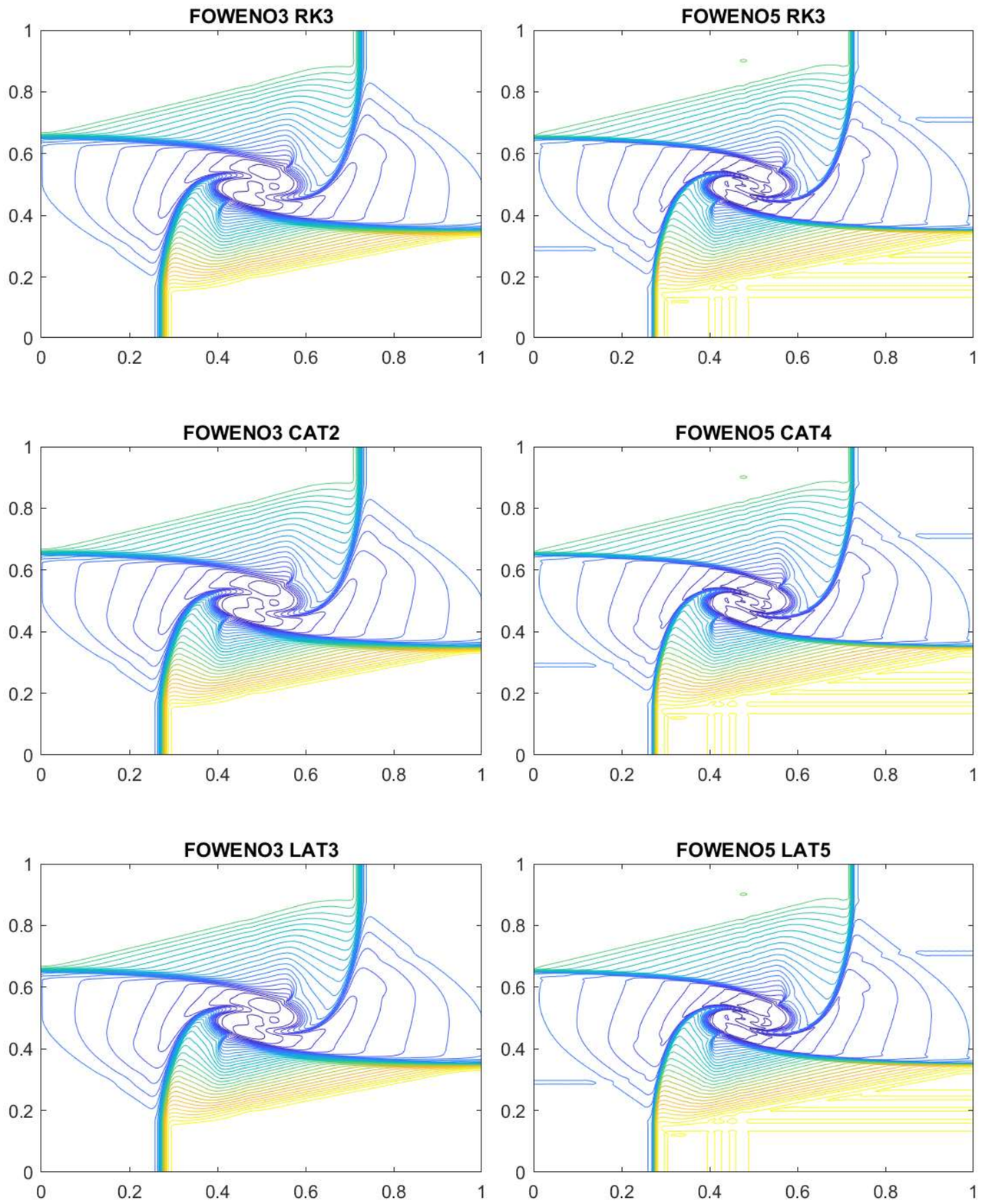}
	\vspace{ -1 cm}
	\caption{Test 9: 2D Euler equations. Lax configuration 6: density computed with  FOWENO-RK, FOWENO-CAT and  FOWENO-LAT.}
    \label{2d_test_6}
	\end{figure}

\begin{figure}[!ht]
	\small
	\setlength{\unitlength}{1mm}
	\centering	
	\includegraphics[width=\textwidth, height=19 cm]{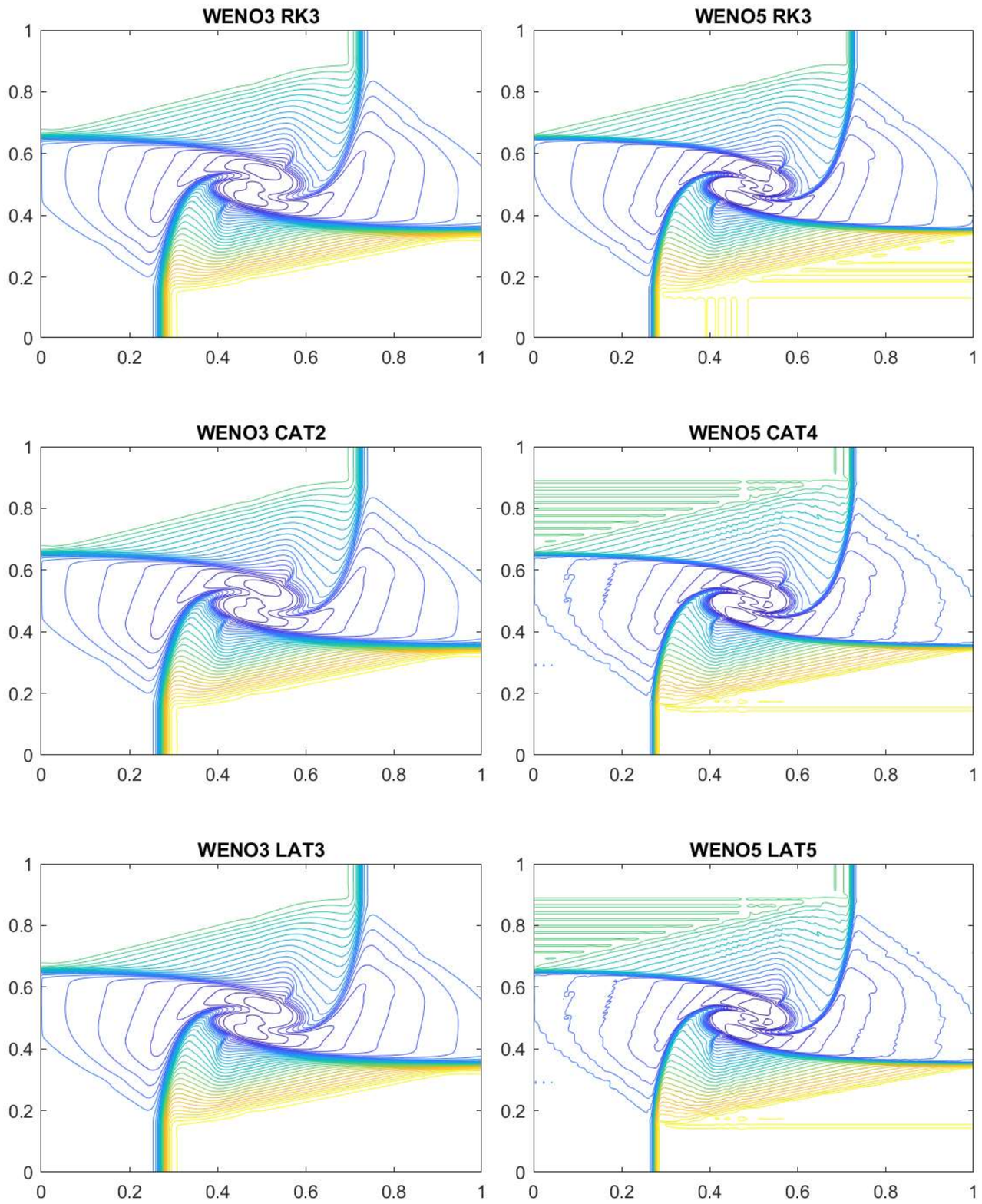}
	\vspace{ -1 cm}
	\caption{Test 9: 2D Euler equations. Lax configuration 6: density computed with  WENO-RK, WENO-CAT and  WENO-LAT.}
    \label{2d_test_6_2}
	\end{figure}

\begin{figure}[!ht]
	\small
	\setlength{\unitlength}{1mm}
	\centering
   \includegraphics[width=\textwidth, height=19 cm]{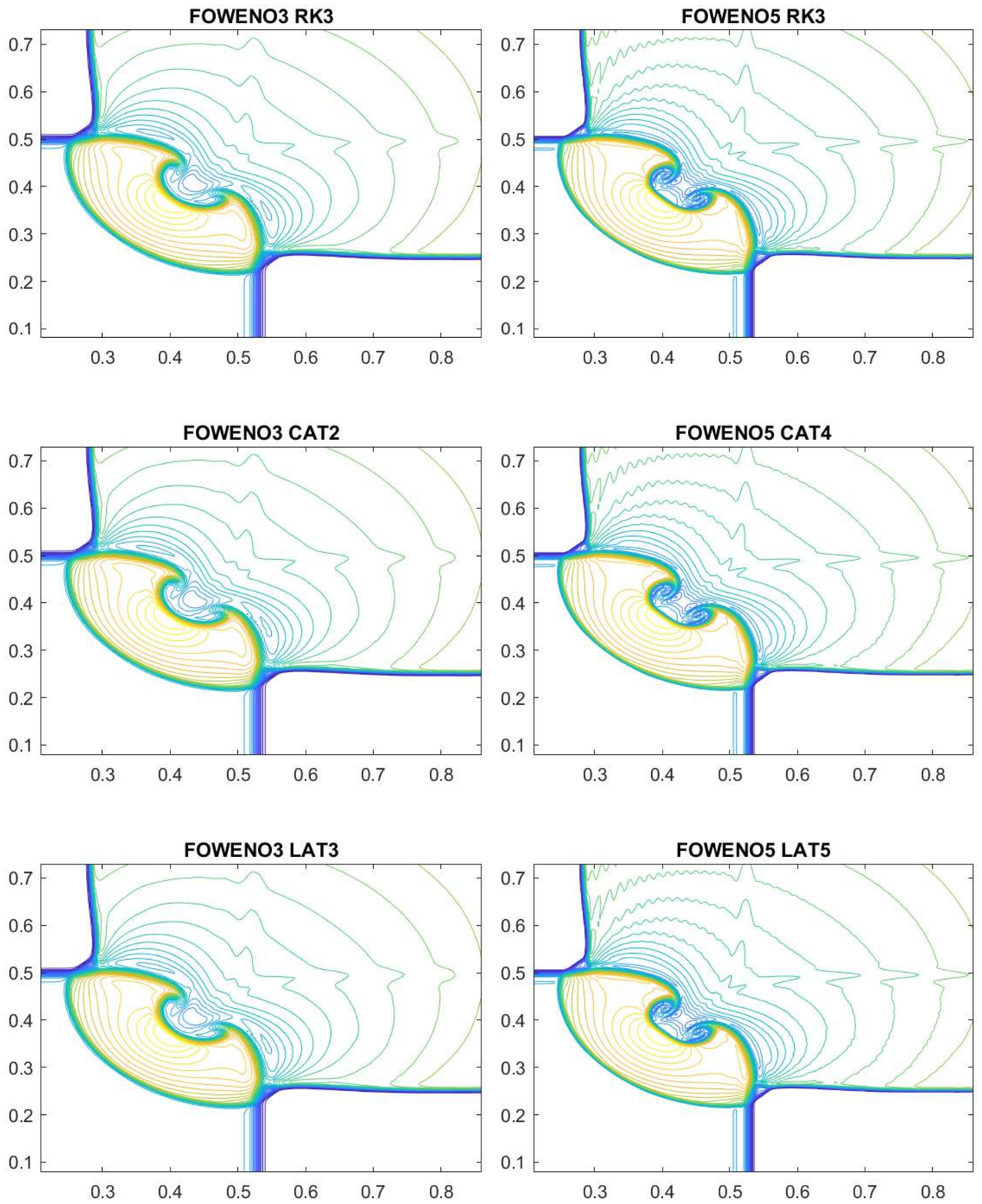}
	\vspace{-1 cm}
	\caption{Test 10: 2D Euler equations. Lax configuration 11: density computed with  FOWENO-RK, FOWENO-CAT and  FOWENO-LAT.}
    \label{2d_test_11}
	\end{figure}

\begin{figure}[!ht]
	\small
	\setlength{\unitlength}{1mm}
	\centering	
    \includegraphics[width=\textwidth, height=19cm]{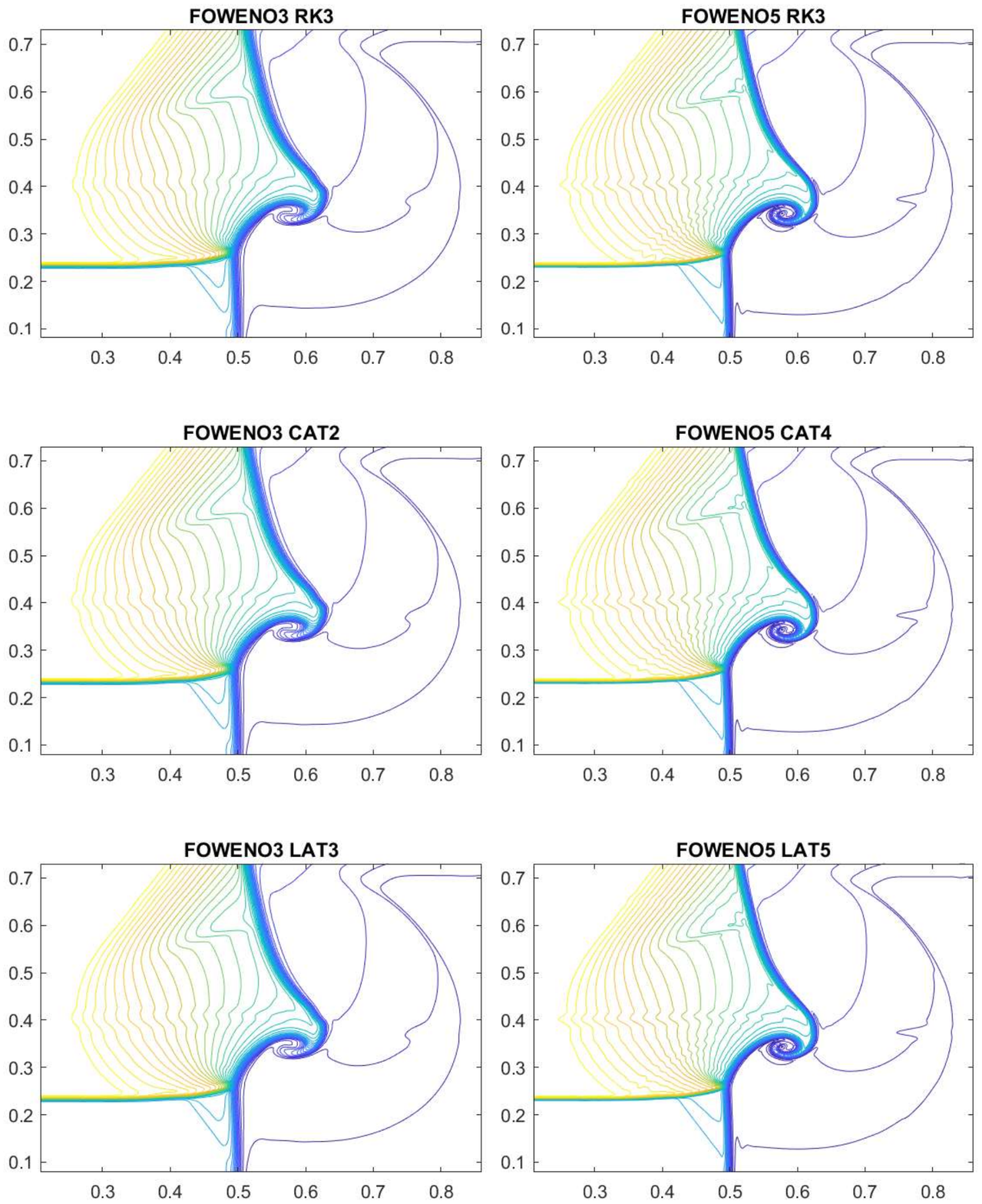}
	\vspace{-1.5 cm}
	\caption{Test 11: 2D Euler equations. Lax configuration 13: density computed with  FOWENO-RK, FOWENO-CAT and  FOWENO-LAT.}
    \label{2d_test_13}
	\end{figure}

\begin{figure}[!ht]
	\small
	\setlength{\unitlength}{1mm}
	\centering	
    \includegraphics[width=\textwidth, height=19cm]{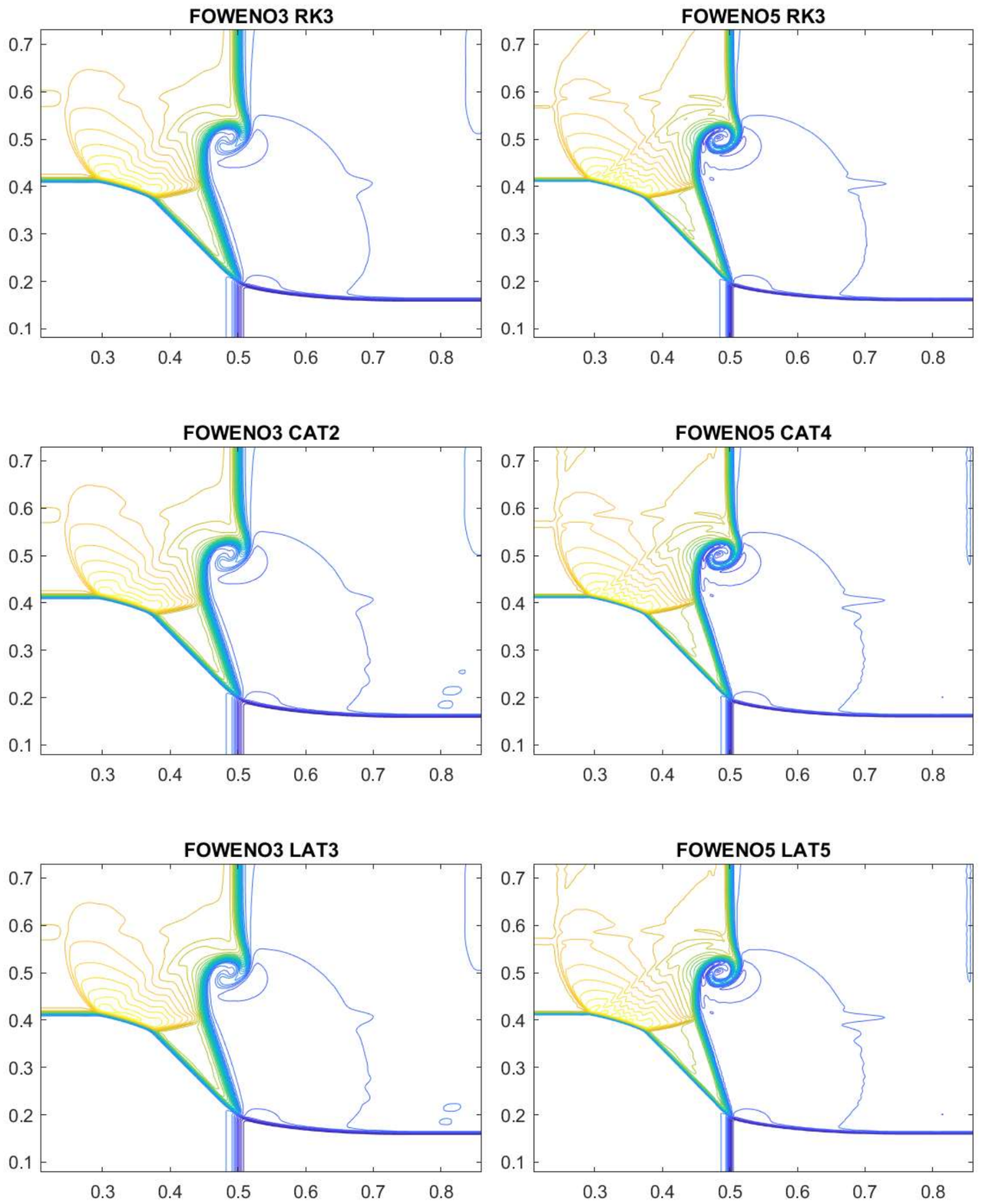}
	\vspace{-1.5 cm}
	\caption{Test 12: 2D Euler equations. Lax configuration 17: density computed with  FOWENO-RK, FOWENO-CAT and  FOWENO-LAT.}
    \label{2d_test_17}
	\end{figure}

\begin{figure}[!ht]
	\small
	\setlength{\unitlength}{1mm}
	\centering
    \includegraphics[width=\textwidth,height=19cm]{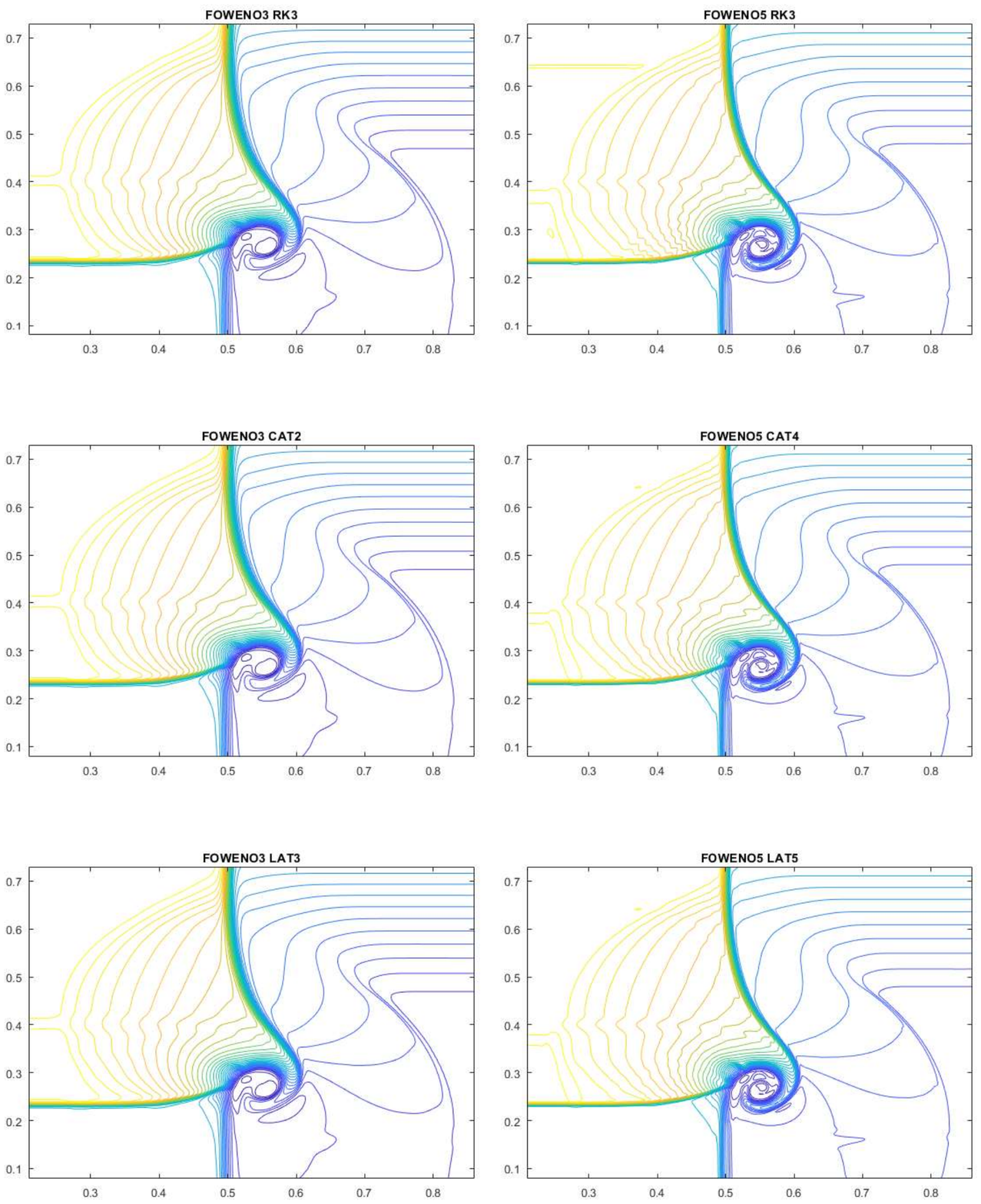}
    \vspace{-1 cm}
	\caption{Test 13: 2D Euler equations. Lax configuration 19: density computed with  FOWENO-RK, FOWENO-CAT and  FOWENO-LAT.}
    \label{2d_test_19}
	\end{figure}

The methods are run in a  $400 \times 400$ point mesh of the computational domain $[0,1] \times [0,1]$ with $CFL=0.475$ and outflow-inflow boundary conditions.  Lax-Friedrichs flux-splitting is used in both WENO and FOWENO implementations. Figures \ref{2d_test_3} to \ref{2d_test_19}  show the numerical densities obtained for the   Lax configurations 3, 6, 11, 13, 17 and 19, respectively. Only the numerical solutions obtained with methods based on FOWENO reconstructions of order 3 or 5 are plotted with the exception of Test 9 for which the solutions given by methods based on WENO reconstructions are also plotted for comparison. Plots are made in Matlab with 25 contour lines.


\begin{table}[htbp]
\begin{center}
\resizebox{\textwidth}{!}{
\begin{tabular}{|llll|rrr|}
\hline
\textbf{Test 8}  &  Configuration 3                &                &  \multicolumn{1}{c}{}  &   \multicolumn{3}{c|}{} \\ 
\hline
 $p_2 = 0.3$     & $\rho_2 = 0.5323$&  $p_1 = 1.5$   & $\rho_1 = 1.5$        &         &        &        \\
 $u_2 = 1.206$   & $v_2=0$          &  $u_1=1$       & $v_1=0$               &         & $\overleftarrow{S_{2,1}}$      &        \\
 $p_3=0.029$     & $\rho_3=0.138$   &  $p_4=0.3$     & $\rho_4=0.5323$       & $\overleftarrow{S_{3,2}}$      &        & $\overleftarrow{S_{4,1}} $    \\
 $u_3=1.206$     & $v_3=1.206$      &   $u_4=0$      & $v_4=1.206$           &          & $\overleftarrow{S_{3,4}}$            &\\
                                      \hline
\end{tabular}
}
\end{center}
\end{table}

\begin{table}[htbp]
\begin{center}
\resizebox{\textwidth}{!}{
\begin{tabular}{|llll|rrr|}
\hline
\textbf{Test 9}  &  Configuration 6                &                &  \multicolumn{1}{c}{}  &   \multicolumn{3}{c|}{} \\ 
\hline
$p_2 = 1$        & $\rho_2 = 2$   &  $p_1 = 1$     & $\rho_1 = 1$   &         &        &        \\
$u_2 = 0.75$     & $v_2=0.5$      &  $u_1= 0.75$   & $v_1=-0.5$     &          & $J^-_{2,1}$      &        \\
$p_3=1$          & $\rho_3=1$     &  $p_4=1$       & $\rho_4=3$     &  $J^+_{3,2}$      &   &     $J^+_{4,1} $    \\
$u_3=-0.75$      & $v_3=0.5$      &  $u_4=-0.75$   & $v_4=-0.5$            &            & $J^-_{3,4}$       & \\      
\hline
\end{tabular}
}
\end{center}
\end{table}

\begin{table}[htbp]
\begin{center}
\resizebox{\textwidth}{!}{
\begin{tabular}{|llll|rrr|}
\hline
\textbf{Test 10}  &  Configuration 11                &                &  \multicolumn{1}{c}{}  &   \multicolumn{3}{c|}{} \\ 
\hline
 $p_2 = 0.4$    & $\rho_2 = 0.5313$    &$p_1 = 1$    & $\rho_1 = 1$                  &          &        &        \\
 $u_2 = 0.8275$ & $v_2=0$     &$u_1= 0.1$    & $v_1=0$                  &          & $\overleftarrow{S_{2,1}}$      &        \\
 $p_3=0.4$ & $\rho_3=0.8$ &  $p_4=0.4$   & $\rho_4=0.5313$               & $J^+_{3,2}$      &        & $\overleftarrow{S_{4,1}}$    \\
 $u_3=0.1$ & $v_3=0$    &  $u_4=0.1$   & $v_4=0.7276$                  &          & $J^+_{3,4}$       &              \\
 \hline
\end{tabular}
}
\end{center}
\end{table}

\begin{table}[htbp]
\begin{center}
\resizebox{\textwidth}{!}{
\begin{tabular}{|llll|rrr|}
\hline
\textbf{Test 11}  &  Configuration 13                &                &  \multicolumn{1}{c}{}  &   \multicolumn{3}{c|}{} \\ 
\hline
$p_2 = 1$    & $\rho_2 = 2$    & $p_1 = 1$    & $\rho_1 = 1$                  &          &        &        \\
$u_2 = 0$ & $v_2=0.3$    & $u_1= 0$    & $v_1=-0.3$                  &          & $J^-_{2,1}$      &        \\
$p_3=0.4$ & $\rho_3=1.0625$ &  $p_4=0.4$   & $\rho_4=0.5313$               &  $\overleftarrow{S_{3,2}}$      &        & $\overleftarrow{S_{4,1}}$    \\
$u_3=0$ & $v_3=0.8145$    &  $u_4=0$   & $v_4=0.4276$                  &          & $J^-_{3,4}$       &   \\
\hline
\end{tabular}
}
\end{center}
\end{table}

\begin{table}[htbp]
\begin{center}
\resizebox{\textwidth}{!}{
\begin{tabular}{|llll|rrr|}
\hline
\textbf{Test 12}  &  Configuration 17                &                &  \multicolumn{1}{c}{}  &   \multicolumn{3}{c|}{} \\ 
\hline
$p_2 = 1$    & $\rho_2 = 2$    & $p_1 = 1$    & $\rho_1 = 1$                   &         &        &        \\
$u_2 = 0$ & $v_2=-0.3$    &  $u_1= 0$    & $v_1=-0.4$                  &          & $J^-_{2,1}$      &        \\
$p_3=0.4$ & $\rho_3=1.0625$ & $p_4=0.4$   & $\rho_4=0.5197$               & $\overleftarrow{S_{3,2}}$      &        & $\overrightarrow{R_{4,1}}$    \\
$u_3=0$ & $v_3=0.2145$    &   $u_4=0$   & $v_4=-1.1259$                  &          & $J^-_{3,4}$       & \\
\hline
\end{tabular}
}
\end{center}
\end{table}

\begin{table}[htbp]
\begin{center}
\resizebox{\textwidth}{!}{
\begin{tabular}{|llll|rrr|}
\hline
\textbf{Test 13}  &  Configuration 19                &                &  \multicolumn{1}{c}{}  &   \multicolumn{3}{c|}{} \\ 
\hline
 $p_2 = 1$    & $\rho_2 = 2$    &  $p_1 = 1$    & $\rho_1 = 1$                  &          &        &        \\
 $u_2 = 0$ & $v_2=-0.3$    &  $u_1= 0$    & $v_1=0.3$                  &        & $J^+_{2,1}$      &        \\
 $p_3=0.4$ & $\rho_3=1.0625$ &  $p_4=0.4$   & $\rho_4=0.5197$               &  $\overleftarrow{S_{3,2}}$      &        & $\overrightarrow{R_{4,1}}$    \\
 $u_3=0$ & $v_3=0.2145$    &   $u_4=0$   & $v_4=-0.4259$                  &         &  $J^-_{3,4}$      &\\
 \hline
 \end{tabular}
}
\end{center}
\end{table}

 In all cases methods based on  third order reconstructions give similar solutions to those provided in \cite{KT2002}, even for FOW3C2 in spite of its lower order of accuracy in time. 
Qualitatively, no significant differences between the results obtained using CAT2 or LAT3 are detected. 
Methods based on fifth order reconstructions are sharper in all cases, as  expected. The quality of the solutions obtained with CAT and LAT are mostly identical again. A comparison between Figures \ref{2d_test_6} and \ref{2d_test_6_2} makes noticeable the improvements provided by FOWENO compared to standard WENO. 

Table \ref{table2d1} shows the CPU time rates for Test 9. Again W3C2 is the cheapest one and its CPU time is takes as the reference. For 3rd order methods, FOW3R3 is the most expensive method. However, for 5th order methods FOW5L5 is the cheapest one and W5C4, the most expensive one. 

\begin{table}[htbp]
\begin{center}
\resizebox{\textwidth}{!}{
\begin{tabular}{|c||c||c||c||c||c|}
\hline
W3R3  & W3C2  & W3L3  & W5R3  & W5C4 & W5L5  \\
  2.5269  & 1.0000    &  1.1228   &  4.7006   & 5.5358   & 3.715   \\
\hline
\hline
FOW3R3  & FOW3C2    &  FOW3L3   &   FOW5R3  &  FOW5C4  &  FOW5L5  \\
2.9967  & 1.2697    &  1.8280    & 4.0197   &  5.1386   & 3.3760     \\
\hline 
\end{tabular}
}
\vspace{2mm}
\caption{CPU time rates for 2D numerical solutions of Test 9.}
\label{table2d1}
\end{center}
\end{table}

\section{Conclusions} \label{sec6}
Several shock-capturing high-order finite difference methods for 1d and 2d systems of conservation laws have been presented and compared in a number of test cases. Two different high-order reconstruction operators have been considered: standard WENO and FOWENO operators. The latter combine the use of fast smooth indicators (that coincide with the original smooth indicators in the third order case) and the computation of optimal weights that allow one to preserve the accuracy of the reconstructions close to critical point regardless of their order. For the best of our knowledge, this is the first time that these two techniques have been combined. 

Concerning the time discretization, two family of methods have been considered: SSPRK methods and Approximate Taylor methods. Moreover, two different implementations of the latter are considered: Lax-Wendroff Approximate Taylor and Compat Approximate Taylor methods. The first one is cheaper but the second one uses smaller stencils and the stability properties are better. 

The numerical tests show that, for third reconstructions, FOWENO is more expensive than WENO  due to the computation of the optimal weights, as it happens for CWENO \cite{CWENO}, M-WENO \cite{HENRICK2005542} and other WENO versions. Nevertheless this extra cost is relatively small and it is compensated by the quality of the solutions close to critical points. For order 5 or bigger, methods based on FOWENO reconstructions give better solutions and are cheaper than those based on standard WENO: the extra cost due to the computation of the optimal weights is compensated by the lower cost required by the computation of the smooth indicators.   

Concerning the time discretization, the following conclusions can be drawn from the numerical tests:
\begin{itemize}

\item CAT2 combined with 3d order reconstructions is a good choice in 1d and 2d: the quality of the solutions is comparable to those obtained with LAT3 or RK3, but with a significantly lower cost. 

\item LAT methods are cheaper for reconstructions of order 7 or bigger in 1d and of order 5 or bigger in 2d, LAT methods. 

\item In some cases, the extra cost of CAT methods  of higher order can be compensated by the fact that bigger values of the CFL parameter can be taken with good results. 

\item For 1d problems, SSPRK3 gives results that are competitive both in quality and computational time. SSPRK4 increases a lot the computational time.

    

\end{itemize}

Approximate Taylor methods are highly parallelizable: future work includes the parallel implementation of these methods in GPU. Another foreseen  extension is the application of Approximate Taylor techniques to obtain high-order well-balanced methods for systems of balance laws.

\section*{Acknowledgements} 
This research has received funding from the European Union’s Horizon 2020 research and innovation program, under the Marie Sklodowska-Curie grant agreement No 642768. It  has been also partially supported by the Spanish Government and FEDER through the Research project RTI2018-096064-B-C21. D. Zor\'io is also supported by Fondecyt Project 3170077.

\bibliographystyle{unsrt}
\bibliography{mybibfile}

\end{document}